\PassOptionsToPackage{tbtags}{amsmath}
\documentclass{ws-m3as}

\usepackage[super]{cite}

\usepackage{mathrsfs}
\usepackage{placeins}
\usepackage{xcolor}
\usepackage[verbose,hypertexnames=false]{hyperref}
\hypersetup{colorlinks=false,allbordercolors=blue,pdfborderstyle={/S/U/W 1}}

\newcommand{\nn}{\nonumber}
\newcommand{\eps}{\varepsilon}
\newcommand{\bT}{\mathbf{T}}
\newcommand{\bB}{\mathbf{B}}
\newcommand{\bU}{\mathbf{U}}
\newcommand{\bx}{\mathbf{x}}

\begin{document}

\markboth{Y. Feng, M. Lin and Y. Wang}{Improved error bounds for Dirac}

\catchline{}{}{}{}{}

\title{Improved error bounds on the Lie-Trotter splitting method for the nonlinear Dirac equation in the massless nonrelativistic regime}

\author{Yue Feng}
\address{School of Mathematics and Statistics, Xi'an Jiaotong University\\
710049 Xi'an, China\\
yue.feng@xjtu.edu.cn}

\author{Meiwen Lin}
\address{School of Mathematics and Statistics, Xi'an Jiaotong University\\
710049 Xi'an, China\\
LMW@stu.xjtu.edu.cn}

\author{Yan Wang}
\address{School of Mathematics and Statistics and Hubei Key Laboratory of Mathematical Sciences,\\
Central China Normal University, 430079 Wuhan, China\\
wang.yan@ccnu.edu.cn}

\maketitle

\begin{abstract}
		We establish improved error bounds for the Lie-Trotter splitting method applied to the nonlinear Dirac equation (NLDE) in the massless nonrelativistic regime characterized by the dimensionless parameter $\eps\in (0, 1]$. In this regime, the solution of the NLDE exhibits highly oscillatory behavior with temporal wavelength $O(\eps)$ and spatial wavelength $O(1)$, together with rapidly propagating outgoing waves traveling at speed $O(1/\eps)$. The NLDE is discretized by the Lie-Trotter splitting method in time and the Fourier pseudospectral method in space, yielding the time-splitting Fourier pseudospectral (TSFP) method. By employing the regularity compensation oscillation (RCO) technique, in which high-frequency modes are controlled through the regularity of the exact solution, while low-frequency modes are analyzed by the phase cancellation and energy method, we derive improved error bounds at $O(\tau)$ for the semi-discrete Lie-Trotter splitting scheme and $O(h^{m-1}+\tau)$ for the fully discrete TSFP method. Finally, numerical experiments are presented to validate the theoretical error bounds and demonstrate their sharpness. In addition, we illustrate the dynamics of the two-dimensional nonlinear Dirac equation with a honeycomb lattice potential.
	\end{abstract}



	\section{Introduction}
\label{sec:intro}

In this paper, we consider the following nonlinear Dirac equation (NLDE) on the torus $\mathbb{T}^d$ ($d=1,2,3 $)  \cite{Dirac1,Dirac2,Thaller}: 
\begin{equation}
\label{eq:Dirac}
i\hbar\partial_t\Psi =  \left(- ic\hbar\sum_{j = 1}^{3}
	\alpha_j\partial_j + mc^2
 \beta\right)\Psi+ e V\Psi+F(\Psi)\Psi, \qquad \bx \in \mathbb{T}^d,
\end{equation}
where $\Psi:= \Psi(t, \bx) = (\Psi_1(t, \bx), \Psi_2(t, \bx), \Psi_3(t, \bx), \Psi_4(t, \bx))^T 
\in \mathbb{C}^4$ is a complex-valued spinor wave function, $t$ is time, $\bx = (x_1, x_2, x_3)^T$ is the spatial coordinate, $\partial_j=\frac{\partial}{\partial_{x_j}}(j=1, \ldots, d)$, $V:=V(\mathbf{x})$  is the real-valued electrical potential. There are many constants including $i=\sqrt{-1}$, $\hbar$ the Planck constant, $m$ the mass, $c$ the speed of light and $e$ the unit charge. We take the nonlinearity as $F(\Psi)=\lambda(\Psi^*\beta \Psi)\beta$ with $\lambda \in \mathbb{R}$ denoting the strength of the nonlinear interaction \cite{LAM} and $\Psi^*=\overline{\Psi}^T$, while $\overline{\Psi}$ denotes the complex conjugate of $\Psi$. In addition, $\beta$ and $\alpha_j \  (j = 1, 2, 3)$ are $4 \times 4$ Dirac matrices of the four-dimensional Clifford algebra
\begin{equation*}
\beta = \begin{pmatrix} I_2 & \mathbf{0} \\ \mathbf{0} & -I_2 \end{pmatrix}, \quad
\alpha_j = \begin{pmatrix} \mathbf{0}  & \sigma_j \\ \sigma_j  &  \mathbf{0}  \end{pmatrix}, \quad j = 1, 2, 3,
\end{equation*}
where $I_n$ is the $n \times n$ identity matrix, and $\sigma_j\ (j = 1, 2, 3)$ are the Pauli matrices defined as
\begin{equation*}
\label{Pauli}
\sigma_1 = \begin{pmatrix} 0 &  1\\ 1 & 0 \end{pmatrix}, \quad
\sigma_2 = \begin{pmatrix} 0 & -i \\ i & 0 \end{pmatrix}, \quad
\sigma_3 = \begin{pmatrix} 1 & 0 \\ 0 &  -1 \end{pmatrix}.
\end{equation*}

The Dirac equation is a fundamental relativistic wave equation that describes the dynamics of elementary spin-1/2 massive particles, such as positrons, electrons and quarks. As a cornerstone of modern particle physics, it provides a successful unification of quantum mechanics and special relativity \cite{ACD}. One of its most remarkable predictions is the existence of antimatter, arising from its negative-energy solutions. This prediction was spectacularly confirmed by the experimental discovery of the positron only a few years after the equation was formulated. In recent decades, the Dirac equation has attracted renewed and widespread interest, driven by both theoretical advances and experimental developments in a variety of fields, including the study of graphene and graphite \cite{NGPNG,NGMJKGDF,SGMH}, intense laser-matter interactions \cite{FLB}, and two- and three-dimensional topological insulators \cite{QZ}.

Introducing the nondimensionalization
\begin{equation}
\label{eq:nd}
    \widetilde{\bx} = \frac{\bx}{x_s},\ \widetilde{t} = \frac{t}{t_s}, \ \widetilde{m} = \frac{m}{m_s}, \ \widetilde{V} = \frac{V}{A_s}, \ \widetilde{\Psi}(\widetilde{t}, \widetilde{\bx}) = \frac{\Psi(t, \bx)}{\Psi_s},
\end{equation}
where $x_s$, $t_s$ and $m_s$ are length unit, time unit and mass unit, respectively. Plugging \eqref{eq:nd} into the NLDE \eqref{eq:Dirac} and taking $\Psi_s = x_s^{-3/2}$ and $A_s = \frac{m_s x_s^2}{et_s^2}$, after some simplification and dropping the   $\widetilde{}$, we obtain the dimensionless nonlinear Dirac equation \cite{HWYY}
\begin{equation}
\label{eq:Dirac_nd}
i\partial_t\Psi=  \left(- \frac{i}{\eps}\sum_{j = 1}^{3}
	\alpha_j\partial_j + \delta\beta\right)\Psi+ V\Psi+F(\Psi)\Psi, \qquad \bx \in \mathbb{T}^d,
\end{equation}
where 
\begin{equation}
    F(\Psi)=\lambda(\Psi^*\beta \Psi)\beta,\qquad \Psi \in \mathbb{C}^4,
\end{equation}
with $\delta = m_0/\eps^2$. Note that $0 < \eps, m_0 \leq 1$ the dimensionless parameters are defined by
\begin{equation}
    \eps := \frac{x_s}{t_s c} = \frac{v_s}{c},\quad m_0 := \frac{m}{m_s}, 
\end{equation}
 where $v_s = x_s/t_s$ is the dimensionless velocity unit, $\eps$ is the ratio between the wave velocity and the speed of light—i.e. it is inversely proportional to the speed of light, and $m_0$ is the ratio between the mass of the particle $m$ and the dimensionless mass unit $m_s$.

In this paper, we mainly consider the massless nonrelativistic regime, i.e., the case where $\delta = O(1) \ (m_0 = O(\eps^2))$, $0 < \eps \ll 1$. 
We focus on the following nonlinear Dirac equation
\begin{equation}
\label{eq:Dirac_nd1}
i\partial_t\Psi=  \left(- \frac{i}{\eps}\sum_{j = 1}^{3}
	\alpha_j\partial_j + \delta\beta\right)\Psi+ V\Psi+F(\Psi)\Psi, \qquad \bx \in \mathbb{T}^3,
\end{equation}
with the initial data
 \begin{equation*}
 	\Psi(0,\bx) =  \Psi_0(\bx), \qquad \bx \in \mathbb{T}^3.
 \end{equation*}
 
More generally, the NLDE in three dimensions (3D) and its dimension reduction to two dimensions (2D) and one dimension (1D) can be expressed in a unified way as 
\begin{equation}
\label{eq:Dirac_nd2}
i\partial_t\Psi=  \left(- \frac{i}{\eps}\sum_{j = 1}^{d}
	\alpha_j\partial_j + \delta\beta\right)\Psi+ V\Psi+F(\Psi)\Psi, \qquad \bx \in \mathbb{T}^d,
\end{equation}
where $d = 1, 2, 3$ indicates the dimension, and $\bx$ is set to be $(x_1, x_2)^T$ in 2D and $x_1$ in 1D. In 1D and 2D, similar to the process in \cite{BCJT}, the nonlinear Dirac equation \eqref{eq:Dirac_nd1} can be decoupled into simplified PDEs with a two-component wave function $\Phi := \Phi(t, \bx) = (\Phi_1(t, \bx), \Phi_2(t, \bx))^T \in \mathbb{C}^2$ satisfying 
\begin{equation}
\label{eq:Dirac_21}
i\partial_t \Phi =  \left(- \frac{i}{\eps}\sum_{j = 1}^{d}
	\sigma_j\partial_j + \delta\sigma_3\right)\Phi + V\Phi+F(\Phi)\Phi, \qquad \bx \in \mathbb{T}^d,
\end{equation}
where 
\begin{equation}
    F(\Phi)=\lambda (\Phi^*\sigma_3 \Phi)\sigma_3,\qquad \Phi \in \mathbb{C}^2.
\end{equation}
  We impose the initial condition
 \begin{equation*}
 	\Phi(0,\bx) =  \Phi_0(\bx), \qquad \bx \in \mathbb{T}^d.
 \end{equation*}
In addition, the NLDE \eqref{eq:Dirac_21} conserves the total mass as
\begin{align}
    \|\Phi(t,\cdot)\|^2:&=\int_{\mathbb{T}^d}|\Phi(t,\bx)|^2d \bx=\int_{\mathbb{T}^d}\sum_{j=1}^2|\Phi_j(t,\bx)|^2d\bx \equiv \|\Phi(0,\cdot)\|^2\nn
    \\ &=\|\Phi_0\|^2,\qquad t \geq 0,
\end{align}
and the energy as
\begin{align}
       E(t):&=\int_{\mathbb{T}^d}\left[-\frac{i}{\eps}\sum_{j=1}^d \Phi^*\sigma_j\partial_{x_j}\Phi+\delta \Phi^*\sigma_3 \Phi + V(\bx)|\Phi|^2+G(\Phi) \right]d\bx,\nn \\
   & \equiv E(0),\qquad  t\geq 0,
\end{align}
 where
\begin{equation}
    G(\Phi)=\frac{\lambda}{2}(\Phi^*\sigma_3\Phi)^2,   \qquad \Phi \in  \mathbb{C}^2.
\end{equation}
Then we introduce the total density $\rho:=\rho(t,\bx)$ as
\begin{equation}
    \rho(t,\bx)=\sum^2_{j=1}\rho_j(t,\bx)=\Phi(t,\bx)^*\Phi(t,\bx), \qquad \bx\in \mathbb{T}^d.
\end{equation}
               
\begin{figure}[!t]
\centering
\includegraphics[width=0.7\textwidth]{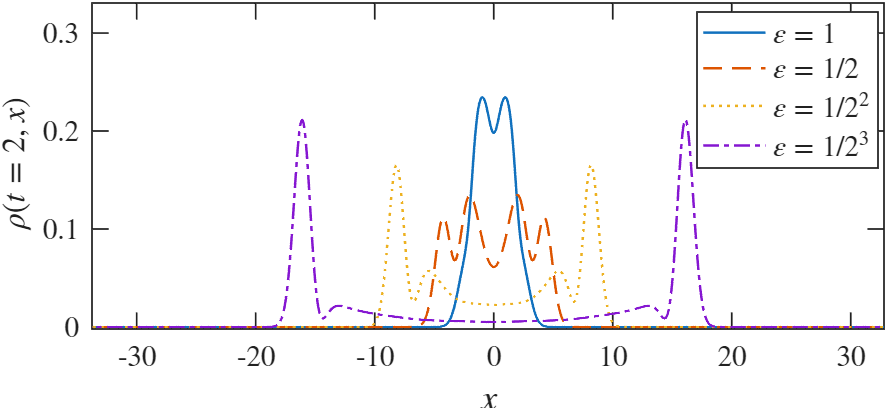}
\caption{Dynamics of total density $\rho$ of the NLDE \eqref{eq:Dirac_21} in 1D for different $\eps$.}
\label{fig:comparison_scalings}
\end{figure}

Over the past decades, various numerical methods have been developed for solving the NLDE with $\eps = 1$. Rigorous error estimates have been established for different numerical methods including the finite difference time domain (FDTD) methods \cite{BCJT, BHM,MY}, the exponential wave integrator Fourier pseudospectral (EWI-FP) method \cite{BCJT,BCJY} and the time-splitting Fourier pseudospectral (TSFP) method \cite{BCY,FLB}. For more details related to the numerical schemes, we refer to \cite{BSG,FLB,Gosse,GSX} and references therein.  When $0<\eps\ll 1$ in the NLDE \eqref{eq:Dirac_21}, i.e., in the massless nonrelativistic regime, the solution exhibits highly oscillatory behavior with temporal wavelength $O(\eps)$ and spatial wavelength $O(1)$. In addition, it generates rapidly propagating outgoing waves with speed $O(1/\eps)$. To illustrate this rapid propagation, Figure~\ref{fig:comparison_scalings} displays the total density $\rho$ of the NLDE at $t=2$ for the one-dimensional case with $\delta=1$, $\lambda=1$ and  
\begin{equation*}
    V(x)=\frac{1}{1+\sin^2(x)}, \qquad \Phi_0(x)=\left(e^{-\frac{x^2}{2}}, 0\right)^T,
\end{equation*}
for different $\eps$. It is well known that such highly oscillatory dynamics pose significant challenges for numerical simulation \cite{BCY2,HWYY}. Standard error estimates for classical numerical methods typically involve negative powers of  $\eps$, implying that their expected convergence rates can only be achieved under suitably restrictive $\eps$-dependent time-step conditions.

Time-splitting methods are particularly attractive for nonlinear Dirac equations, since they allow the linear Dirac operator and the nonlinear or potential terms to be treated separately. When combined with Fourier pseudospectral discretization in space, they yield efficient schemes for periodic problems. However, the standard error analysis for the first-order Lie-Trotter splitting method in the massless nonrelativistic regime still gives the error bound at $O(\tau^2/\eps)$ \cite{HWYY}, which does not provide a satisfactory uniform error bound as $\eps\to 0$. 

For the nonlinear Dirac equation, the RCO technique was recently introduced to establish the improved uniform error bounds for the TSFP method in the long-time regime \cite{BCF2}. Motivated by this development, the aim of this paper is to establish the improved error bounds for the Lie-Trotter splitting Fourier pseudospectral method for the NLDE \eqref{eq:Dirac} in the massless nonrelativistic regime. Following the RCO framework, we introduce a frequency cut-off parameter $\tau_0$. The high-frequency modes with frequencies larger than $(1/\tau_0)$ are controlled by exploiting the regularity of the exact solution, while the low-frequency modes with frequencies not exceeding $ (1/\tau_0)$ are analyzed through a combination of phase-cancellation arguments and energy method.

The remainder of this paper is organized as follows. In Section 2, we discretize the NLDE \eqref{eq:Dirac_21} in time using the Lie-Trotter splitting method to obtain the semidiscrete scheme and then combine it with the Fourier pseudospectral discretization in space to derive the fully discrete scheme. In Section 3, we establish improved error bounds for both the semidiscrete and fully discrete schemes. Numerical results are presented in Section 4. Finally, some conclusions are drawn in Section 5. Throughout this paper, the notation $A\lesssim B$ signifies that there exists a generic constant $C>0$, independent of the mesh size $h$, the time step $\tau$, and the parameter $\eps$ such that $\vert A\vert \leq CB$.

\section{Lie-Trotter Splitting Fourier Pseudospectral Method} \label{sec:TSFP}
	For simplicity, we present the numerical method and the corresponding error estimates for the NLDE \eqref{eq:Dirac_21} in 1D, i.e., $d=1$. The methods and results can be easily generalized to \eqref{eq:Dirac_21} in 2D and 3D with minor modifications. The NLDE \eqref{eq:Dirac_21} in 1D on the bounded domain $\Omega = (a, b)$ with periodic boundary conditions is given by
	\begin{align}
		\label{eq:Dirac_1D}
		  & i\partial_t \Phi = \left(- \frac{i}{\eps} \sigma_1\partial_x + \delta \sigma_3  \right)\Phi + V\Phi + \lambda (\Phi^*\sigma_3 \Phi)\sigma_3\Phi,\qquad  x\in\Omega, \quad t>0, \\
		\label{eq:ib}
	  &  \Phi(t,a)=\Phi(t,b),  \quad t\geq0 ; \quad \Phi(0,x)=\Phi_0(x), \qquad x\in\overline{\Omega},
	\end{align}
    where $\Phi:= \Phi(t, x)$ and $\Phi_0(a) = \Phi_0(b)$.
	
	Introducing the following operators
\begin{equation}\label{eq:operators}
		\bT := - \frac{i}{\eps} \sigma_1\partial_x + \delta\sigma_3 , \qquad \bB(\Phi) := V I_2 + \lambda \bigl(\Phi^{*}\sigma_3 \Phi\bigr)\sigma_3,
	\end{equation}
	we can express \eqref{eq:Dirac_1D} as
	\begin{equation}
		i\partial_t \Phi(t, x) =\bT \Phi(t,x)+\bB(\Phi(t,x))\Phi(t,x) , \qquad x\in\Omega,\ t>0.
		\label{eq:Dirac1d}
	\end{equation}
	
	Choose the time step size $\tau := \Delta t>0$ and time steps $t_n=n\tau$ for $n = 0, 1, 2, \ldots$. Let $\mathcal{S}_{\rm exact}^{\tau}$ denote the exact flow of \eqref{eq:Dirac1d} in one time step, then the exact solution is
    \begin{equation}
        \Phi(t_{n+1})=\mathcal{S}_{\rm exact}^{\tau}(\Phi(t_n)), \qquad n=0, 1, 2, \ldots.
    \end{equation}

    To approximate the exact flow $\mathcal{S}_{\rm exact}^\tau$, we apply the first-order Lie-Trotter splitting \cite{Trotter}. Denote by $\mathcal{S}T^\tau$ and $\mathcal{S}_B^\tau$ the exact flows associated with the following two subproblems:
    \begin{align*}
       \varphi_T^\tau: \ i\partial_t \Phi &= \bT \Phi,  \\
       \varphi_B^\tau:\  i\partial_t \Phi &=\bB(\Phi)\Phi,
    \end{align*}
 then the Lie-Trotter splitting flow is defined by 
    \begin{equation}
        \mathcal{S}_{\tau} := \mathcal{S}_T^\tau\circ\mathcal{S}_B^\tau,
    \end{equation}
which means that, in each time step, the nonlinear subproblem is solved first, followed by the linear Dirac subproblem.
    
    Denoting by $\Phi^{[n]}$ the approximations of $\Phi(t_n)$ and taking the initial condition $\Phi ^{[0]}(x):=\Phi_0(x)$ for $x\in\Omega$, the semi-discretization of the NLDE \eqref{eq:Dirac_1D} via Lie-Trotter splitting is given by
	\begin{equation}\label{eq:S1}
		\Phi ^{[n+1]}= \mathcal{S}_{\tau}(\Phi^{[n]}) =\mathcal{S}_T^\tau\left(\mathcal{S}_B^\tau(\Phi^{[n]})\right), \qquad n=0, 1, 2, \ldots.
	\end{equation}

	\subsection{Full-discretization}
	In this subsection, we present the fully discrete implementation of the Lie-Trotter splitting method in time combined with the Fourier pseudospectral method in space.

	For an integer $m \geq 0$, we denote by $H^m_{\rm per}(\Omega)=\{\ u \in H^m(\Omega) | \partial^k_x u(a) = \partial^k_x u(b), k = 0, 1, \ldots, m-1 \}$ with the $H^m$-norm on $H^m_{\rm per}(\Omega)$ given by
	\begin{equation*}
		\|u\|_{H^m}=\left(\sum_{l \in \mathbb{Z}} \left(1+\mu_l^2\right)^m \left|\widehat{u}_l\right|^2\right)^{1/2}, \quad \text{for} \quad  u(x)=\sum_{l \in \mathbb{Z}} \widehat{u}_l e^{i\mu_l(x-a)}, \quad \mu_l=\frac{2\pi l}{b-a},		
\end{equation*}
    where $\widehat{u}_l \, (l \in \mathbb{Z})$ are the Fourier coefficients of the function $u(x)$  \cite{ST}. 

    Let $M$ be an even positive integer and choose the spatial mesh size $h = (b-a)/M$, then the grid points are given as
	\begin{equation*}
		x_j = a + jh,  \quad j=0, 1, \ldots, M.
	\end{equation*}
We denote the sets  $\mathcal{T}_M $, $X_M$, $Y_M$, $Z_M$ as
	\begin{align*}
		& X_M  = \left\{U = (U_0, U_1, \ldots, U_M)^T~|~ U_j\in\mathbb{C}^2, \, j=0, 1, \ldots, M, \, U_0 = U_M\right\},\\
		& Y_M = Z_M\times Z_M, \quad Z_M = {\rm{span}}\{\phi_l(x)=e^{i\mu_l(x-a)}, \, l\in\mathcal{T}_M\},
	\end{align*}
    where the index set $\mathcal{T}_M= \{l~|~l=-M/2, -M/2+1, \ldots, M/2-1\}$.
The projection operator $P_M: \left(L^2(\Omega)\right)^2\rightarrow Y_M$ is defined as
	\begin{equation*}
		\left(P_MU\right)(x) := \sum_{l\in\mathcal{T}_M}\widehat{U}_le^{i\mu_l(x-a)}, \quad U(x)\in \left(L^2(\Omega)\right)^2,
	\end{equation*}
	where
	\begin{equation}
		\widehat{U}_l = \frac{1}{b-a}\int_a^bU(x)e^{-i\mu_l(x-a)}dx, \quad l\in\mathcal{T}_M.
	\end{equation}
	Define the space $\left(C_{\rm per}(\overline{\Omega})\right)^2 = \left\{U \in \left(C(\overline{\Omega})\right)^2~|~U(a) = U(b) \right\}$ and the interpolation operator $I_M : \left(C_{\rm per}(\overline{\Omega})\right)^2\rightarrow Y_M$ or $I_M: X_M\rightarrow Y_M$ as
	\begin{equation*}
		\left(I_MU\right)(x) := \sum_{l\in\mathcal{T}_M}\widetilde{U}_le^{i\mu_l(x-a)}, \qquad  U(x)\in \left(C_{\rm per}(\overline{\Omega})\right)^2\quad {\rm{or}}\quad U\in X_M,
	\end{equation*}
	where
	\begin{equation}
		\widetilde{U}_l = \frac{1}{M}\sum_{j=0}^{M-1}U_je^{-2ijl\pi/M}, \quad l\in\mathcal{T}_M,
	\end{equation}
	with $U_j = U(x_j)$ for a function $U(x)$.
	
	Let $\Phi^n = \left(\Phi_0^n, \Phi_1^n, \ldots, \Phi_M^n\right)^T \in X_M$ denote the solution vector in $t=t_n$ with $\Phi_j^n$ the numerical approximation of $\Phi(t_n, x_j)$. Taking the initial value $\Phi_j^0=\Phi_0(x_j)$ for $j=0, 1, \ldots, M$, then the Lie-Trotter Fourier pseudospectral (LTFP) method is given to discretize the NLDE \eqref{eq:Dirac_1D}.
	\begin{equation}
		\begin{split}
				&  \Phi^{(1)}_j = e^{-i\tau \Lambda_j^n}\Phi^n_j,\\
                & \Phi^{n+1}_j = \sum_{l \in \mathcal{T}_M} e^{-i\tau \Gamma_l}\widetilde{(\Phi^{(1)})}_l\; e^{i\mu_l(x_j-a)} =  \sum_{l \in \mathcal{T}_M} Q_l\; e^{-i\tau D_l}(Q_l)^T \widetilde{(\Phi^{(1)})}_l\; e^{\frac{2ijl\pi}{M}},
        \end{split}
		\label{eq:TSFP}
	\end{equation}
for $n=0,1,\ldots$, where $\Gamma_l = \frac{\mu_l}{\eps}\sigma_1 + \delta\sigma_3 = Q_l D_l(Q_l)^T$ with $\omega_l =\sqrt{\delta^2+\frac{\mu^2_l}{\eps^2}}$, $(Q_l)^T$ is the transpose of $Q_l$ and
	\begin{equation*}
		\Gamma_l = \begin{pmatrix} \delta &\   \frac{\mu_l}{\eps} \\ \frac{\mu_l}{\eps} &\  -\delta \end{pmatrix}, \quad
		Q_l = \begin{pmatrix} \frac{\omega_l+\delta}{\sqrt{2\omega_l(\omega_l+\delta)}}  & \  -\frac{\frac{\mu_l}{\eps}}{\sqrt{2\omega_l(\omega_l+\delta)}} \\ \frac{\frac{\mu_l}{\eps}}{\sqrt{2\omega_l(\omega_l+\delta)}} &\    \frac{\omega_l+\delta}{\sqrt{2\omega_l(\omega_l+\delta)}} \end{pmatrix}, \quad
		D_l = \begin{pmatrix} \omega_l &\ 0 \\ 0 &\  -\omega_l \end{pmatrix},
	\end{equation*}
	and $\Lambda_j = \textrm{diag}(\Lambda_{j, +},\Lambda_{j, -})$ with $\Lambda^n_{j,\pm}=V(x_j)\pm \lambda(\Phi_j^n)^*\sigma_3 \Phi_j^n$.

\section{Improved Error Bounds} \label{sec:unibds}

In this section, we rigorously establish the improved error bounds for the Lie-Trotter splitting method applied to the NLDE in the massless nonrelativistic regime over a fixed time interval $[0, T]$, where $T>0$.

        \subsection{Main results}
 Let $0<T<T_{\max}<\infty$ with $T_{\max}$ being the maximal existence time  (uniform in $\varepsilon\in(0,1]$) for the solutions to the NLDE \eqref{eq:Dirac_1D}. We make the following assumptions on the electric potential $V:= V(x)$ and the exact solution $\Phi:=\Phi(t, x)$:
    \begin{align*}
        & \text{(A)}\quad \|V\|_{W_{\rm per}^{m, \infty}} \lesssim 1, \qquad m \geq 3,
   \\
       & \text{(B)}\quad
        \|\Phi\|_{L^\infty([0, T]; (H_{\rm per}^{m})^2)} \lesssim 1,\qquad
      \quad    m\geq 3. 
    \end{align*}
     
    \begin{theorem}\label{thm:S2}
         Under assumptions (A) and (B), let $0<\tau_0<1$ be sufficiently small and independent of $\eps$. Then, for any fixed constant $\alpha \in (0,1)$, whenever  $0<\tau \le \eps\alpha \frac{ \pi (b-a)\tau_0}{2\sqrt{\delta^2 \eps^2 \tau_0^2 (b-a)^2+4\pi^2(1+\tau_0)^2}}$, the following improved error bound holds for any $\eps \in (0,1]$
    \begin{equation}\label{eq:semi-cutoff}
          \|\Phi(t_n,x)-\Phi^{[n]}\|_{H^1} \lesssim \tau+\tau_0^{m-1},
           \qquad  0 \le n \le \frac{T}{\tau}.
    \end{equation}
    The constant in the error estimate is independent of $\eps$, but the admissible time step is subject to the coupling condition stated above.
    
     In particular, if the exact solution is sufficiently smooth, e.g., $\Phi(t,x)\in (H_{per}^\infty(\Omega))^2$, then the last term $\tau_0^{m-1}$ decays exponentially with respect to $m$ and is negligible in practice for sufficiently small $\tau_0$. Consequently, the improved error bound would become
    
     \begin{equation}\label{Hthm}
        \|\Phi(t_n,x)-\Phi^{[n]}\|_{H^1} \lesssim \tau,
         \qquad  0 \le n \le \frac{T}{\tau}.
     \end{equation}

\end{theorem}

Let $\Phi^n$ denote the numerical approximation obtained by the LTFP method \eqref{eq:TSFP}. Then, the following improved error bounds hold for the full discretization.

 \begin{theorem}\label{thm:S22}
         Under assumptions (A) and (B), sufficiently small constants $h_0>0$ and $0<\tau_0<1$ can be chosen  independently of $\eps$. Then, for any fixed constant $\alpha \in (0,1)$, whenever $0<h\leq h_0$ and  $0<\tau \le \eps\alpha \frac{ \pi (b-a)\tau_0}{2\sqrt{\delta^2 \eps^2 \tau_0^2 (b-a)^2+4\pi^2(1+\tau_0)^2}}$, the following improved error bound holds for any $\eps \in (0,1]$
    \begin{equation}
          \|\Phi(t_n,x)-I_M \Phi^{n}\|_{H^1} \lesssim h^{m-1}+\tau+\tau_0^{m-1},
           \qquad  0 \le n \le \frac{T}{\tau}.
    \end{equation}
The constant in the error estimate is independent of $\eps$, but the admissible time step is subject to the coupling condition stated above.
    
     Similarly, if the exact solution is sufficiently smooth, e.g. $\Phi(t,x)\in (H_{\rm per}^\infty(\Omega))^2$, the improved error bound would become
    
     \begin{equation}\label{Hthm2}
        \|\Phi(t_n,x)-I_M \Phi^{n}\|_{H^1} \lesssim h^{m-1}+\tau,
         \qquad  0 \le n \le \frac{T}{\tau}.
     \end{equation}
The fully discrete estimate follows by combining the time-semidiscrete error bound with a Fourier pseudospectral comparison in $H^1$. We omit details of the spatial comparison.

\end{theorem}



\begin{remark}
   In Theorems~\ref{thm:S2} and \ref{thm:S22}, as well as in the other results of this paper concerning the one-dimensional problem, we establish the error bounds for $\Phi(x,t)$ in the $H^1$-norm due to the fact that $H^r$ is an algebra for $r>d/2$. In 2D and 3D cases, analogous estimates would instead need to be established in the $H^2$-norm(2D case in the sense of \eqref{eq:Dirac_21}, and 3D case in the sense of the four-component NLDE given in \cite{BCJY} ), which require higher regularity assumptions on the exact solution.
\end{remark}

    \subsection{Proof of Theorem 3.1}
    For  simplicity of notation, we shall write $\Phi(t_n):=\Phi(t_n,x)$ for
$x\in \Omega$ in the proof. Denote the exact flow of the solution 
$\Phi(t_n)\to \Phi(t_{n+1})$ as
\begin{equation}
\Phi(t_{n+1})=\mathcal{S}_{\rm exact}^{\tau}\bigl(\Phi(t_n)\bigr),
\qquad
0\le n\le \frac{T}{\tau}-1.
\end{equation}

We begin by analyzing the local truncation error associated with a single time step of the numerical scheme \eqref{eq:S1}. A Taylor expansion then yields the following expression
\begin{equation*}
   \mathcal{S}_{\tau}\bigl(\Phi(t_n)\bigr)= e^{-i\tau \bT}
\left(
\Phi(t_n)-i\tau B^n \Phi(t_n)
-\tau^2 (B^n)^2 \int_0^1 (1-s)e^{-is\tau B^n}\Phi(t_n)\,ds
\right).
\end{equation*}
On the other hand, by applying Duhamel's principle,  we obtain the following representation of the exact solution $\Phi(t_n+s)$ 
\begin{equation*}
    \Phi(t_n+s)
=
e^{-is\bT}\Phi(t_n)
-
i\int_0^s e^{-i(s-w)\bT}B^n(w)\Phi(t_n+w)\,dw,
\end{equation*}
with $B^n(w) :=V I_2+F(\Phi(t_n+w))$ and $B^n:=B^n(0)$.

We next establish the following estimates for the local truncation error.
\begin{lemma}\label{lem:local_error}
    For any \  $0<\eps\le 1$, the local truncation error of the time-discrete  Lie-Trotter splitting method \eqref{eq:S1} can be decomposed as follows
    \begin{align*}
        \mathscr{E}^n
     :&= \mathcal{S}_\tau\bigl(\Phi(t_n)\bigr)-\Phi(t_{n+1}) \notag\\
     &= \mathscr{F}\bigl(\Phi(t_n)\bigr)
     + \mathscr{R}_1^n + \mathscr{R}_2^n + \mathscr{R}_3^n,    
    \end{align*}    
       where
       \begin{align*}
        \mathscr{F}\bigl(\Phi(t_n)\bigr)=\ &i\bigg( \int_0^\tau e^{-i(\tau-s)\bT}\ \Big(V I_2+F(e^{-is\bT }\Phi(t_n))\Big)\ e^{-is\bT}\Phi(t_n)\ ds 
        \\&-\tau \ e^{-i\tau \bT}\ \Big(V I_2+    F(\Phi(t_n))\Big) \Phi(t_n)  \bigg ),\\
        \mathscr{R}_1^n= \ & i \int_0^\tau e^{-i(\tau-s)\bT}\ \Big(F(\Phi(t_n+s))-F(e^{-is\bT }\Phi(t_n))\Big) e^{-is\bT}\Phi(t_n)\ ds  ,\\
        \mathscr{R}_2^n= \ & i \int_0^\tau e^{-i(\tau-s)\bT}\ \Big(V I_2+F(\Phi(t_n+s))\Big)\ \Big(\Phi(t_n+s)-e^{-is\bT}\Phi(t_n)\Big)\ ds ,\\
        \mathscr{R}_3^n=\  &-e^{-i\tau \bT} \tau^2\ (B^n)^2 \int_0^1(1-s)e^{-is\tau B^n}\Phi(t_n)\ ds.
       \end{align*} 
       
       Under assumptions (A) and (B), the following error bounds hold
       \begin{equation}\label{eq:LE}
           \| \mathscr{F}\bigl(\Phi(t_n)\bigr)\|_{H^1} \lesssim  \frac{\tau^2 }{\eps},
           \quad 
           \|\mathscr{R}_1^n\|_{H^1} \lesssim  \tau^2,
             \quad
           \|\mathscr{R}_2^n\|_{H^1} \lesssim \tau^2,
           \quad
        \|\mathscr{R}_3^n\|_{H^1} \lesssim \tau^2.
       \end{equation}   
\end{lemma}

      \begin{proof}
      Under assumptions (A) and (B), we estimate each term in the local truncation error individually. Using the corresponding quadrature rules, we obtain 
      \begin{align*}
          \mathscr{F}\bigl(\Phi(t_n)\bigr)=\ &i\bigg( \int_0^\tau e^{-i(\tau-s)\bT}\ \Big(V I_2+F(e^{-is\bT }\Phi(t_n))\Big)\ e^{-is\bT}\Phi(t_n)\ ds 
        \\&-\tau \ e^{-i\tau \bT}\ \Big(V I_2+    F(\Phi(t_n))\Big) \Phi(t_n)   \bigg )
        \\=\ &i\int_0^\tau \int_0^s \partial_w G(w)dwds,
      \end{align*}
         where $G(w)=e^{-i(\tau-w)\bT}\ \Big(V I_2+F(e^{-iw\bT }\Phi(t_n))\Big)\ e^{-iw\bT}\Phi(t_n)$. In view of $\bT=O(\frac{1}{\eps})$,  assumptions (A) and (B) imply that
        \begin{equation}
       \left \| \mathscr{F}\big(\Phi(t_n)\big)\right\|_{H^1}=\left\| \int_0^\tau \int_0^s \partial_w G(w)dwds\right\|_{H^1} \lesssim  \frac{\tau^2 }{\eps}.
        \end{equation}

By Duhamel’s principle, we obtain
\begin{equation*}
     \Phi(t_n+s)-e^{-is\bT}\Phi(t_n)
=-i\int_0^s e^{-i(s-w)\bT}B^n(w)\Phi(t_n+w)\,dw.
\end{equation*}
With a direct computation and the definition of the nonlinearity $F(\Phi(t_n))$, we have
\ $ \|\mathscr{R}_1^n\|_{H^1} \lesssim  \tau^2$,   $ \|\mathscr{R}_2^n\|_{H^1} \lesssim \tau^2$,
and 
\begin{equation*}
    \| \mathscr{R}_3^n\|_{H^1}=\left\|-e^{-i\tau \bT} \tau^2\ (B^n)^2 \int_0^1(1-s)e^{-is\tau B^n}\Phi(t_n)\ ds\right\|_{H^1} \lesssim \tau^2 \ .
\end{equation*}
   
  Combining the above estimates completes the proof of
Lemma~\ref{lem:local_error}. The lemma gives the one-step
bound $\|\mathscr{E}^n\|_{H^1}\lesssim \tau^2/\eps$, which
implies $\|\mathscr{E}^n\|_{H^1}\lesssim\tau$ under the
time-step restriction $\tau\lesssim\eps$. We use this
estimate in the global error analysis below.
      \end{proof}

Having established the necessary estimates for the local truncation error, we are now ready to prove Theorem~\ref{thm:S2}. With the error function defined by
\begin{equation}
    \mathbf{e}^{[n]} :=\mathbf{e} ^{[n]}(x) = \Phi^{[n]} - \Phi(t_n), \qquad n = 0,1,\ldots,
\end{equation}
the corresponding error equation takes the form
\begin{equation}\label{eq:Dirac_GE}
    \mathbf{e}^{[n+1]}=\mathcal{S}_\tau\bigl(\Phi^{[n]}\bigr)-
    \mathcal{S}_\tau\bigl(\Phi(t_n)\bigr)+\mathscr{E}^n=e^{-i\tau \bT} \mathbf{e}^{[n]}+W^n+\mathscr{E}^n,\qquad n \ge 0,
\end{equation}
where $W^n:=W^n(x)$ is defined by
\begin{equation*}
    \begin{aligned}
        W^n(x)&=e^{-i\tau\bT}\biggl[\left(e^{-i\tau\big( V I_2+F\left(\Phi^{[n]}\right)\big)}- I_2\right)\Phi^{[n]}\\
        &\qquad-\left(e^{-i\tau\big( V I_2+F(\Phi(t_n))\big)}- I_2\right)\Phi(t_n)\biggr].
    \end{aligned}
\end{equation*}
Based on \eqref{eq:Dirac_GE}, a straightforward induction yields the following representation of the error equation
\begin{equation}\label{EE}
    \mathbf{e}^{[n+1]}
= e^{-i(n+1)\tau \bT} \mathbf{e}^{[0]}
+ \sum_{k=0}^{n} e^{-i(n-k)\tau \bT}\left(W^k + \mathscr{E}^k\right),
\qquad 0 \leq n \leq \frac{T}{\tau} - 1.
\end{equation}
Since $\Phi^{[0]}$ coincides with the exact initial value $\Phi(0)$, the initial error vanishes, i.e., $e^{[0]}=0$. The error equation can be rewritten as follows
\begin{equation}\label{EEF}
    \mathbf{e}^{[n+1]}
= \sum_{k=0}^{n} e^{-i(n-k)\tau \bT}\left(W^k + \mathscr{E}^k\right),
\qquad 0 \leq n \leq \frac{T}{\tau} - 1.
\end{equation}
To estimate $\|W^n(x)\|_{H^1}$, it is necessary first to prove that $\|\Phi^{[n]}\|_{H^1}$ remains uniformly bounded. 
We proceed by mathematical induction. For $n=0$, the estimate follows immediately from assumption (B). Indeed,
\begin{equation*}
    \left\|\Phi^0\right\|_{H^1}=\|\Phi(0,x)\|_{H^1}\lesssim 1.
\end{equation*}

Let $m$ be a fixed nonnegative integer satisfying $0\leq m\leq \frac{T}{\tau}-1$. Suppose that for every integer $n$ with $0\leq n\leq m$,
\begin{equation}
    \left\|\Phi^{[n]}\right\|_{H^1}\lesssim 1.
\end{equation}
Then, under this assumption, $W^n$ satisfies
\begin{equation*}
    	\|W^n\|_{H^1}
	\lesssim
	\tau \left\|e^{[n]}\right\|_{H^1},   \qquad 0\leq n\leq m.
\end{equation*}

On the other hand, it follows from the local truncation error estimate and the time-step restriction $\tau\lesssim\eps$ that
\begin{equation*}
    	\|\mathscr{E}^n\|_{H^1}
	\lesssim
	\frac{\tau^2}{\varepsilon}
	\lesssim
	\tau.
\end{equation*}
Substituting the above estimates for $W^n$ and $\mathscr{E}^n$ into \eqref{EE} produces the following estimate
\begin{equation*}
    	\left\|e^{[n+1]}\right\|_{H^1} \lesssim T+\tau \sum_{k=0}^{n}\left\|e^{[k]}\right\|_{H^1}.
\end{equation*}

Applying the discrete Gronwall inequality to the preceding estimate yields
\begin{equation*}
    \left\|e^{[n+1]}\right\|_{H^1}\lesssim 1,  \qquad 0\leq n\leq m+1.
\end{equation*}
 Therefore,
\begin{equation*}
    	\left\|\Phi^{[n+1]}\right\|_{H^1}
	\leq
	\left\|e^{[n+1]}\right\|_{H^1}
	+
	\|\Phi(t_{n+1})\|_{H^1}
	\lesssim 1.
\end{equation*}
Thus, the induction argument is closed and, consequently,
\begin{equation}\label{bound}
    \left\|\Phi^{[n]}\right\|_{H^1}\lesssim 1,     \qquad 0\leq n\leq \frac{T}{\tau}.
\end{equation}
It further follows from the bound \eqref{bound} that
\begin{equation}\label{eq:Wn}
    \left\| W^n(x) \right\|_{H^1} \lesssim  \tau \left\| e^{[n]} \right\|_{H^1}.
\end{equation}

Combining \eqref{eq:LE} and \eqref{eq:Wn}, we have the estimates for $0\leq n \leq \frac{T}{\tau}-1$,
\begin{equation}\label{eq:erfinial}
    \left\| \mathbf{e}^{[n+1]} \right\|_{H^1}
\lesssim
 \tau
+ \tau \sum_{k=0}^{n} \left\|\mathbf{e} ^{[k]} \right\|_{H^1}
+ \left\| \sum_{k=0}^{n} e^{-i(n-k)\tau \bT}\mathscr{F}\bigl(\Phi(t_k)\bigr) \right\|_{H^1}.
\end{equation}
To derive the improved error bounds, we  employ the {\textbf{regularity compensation oscillation}} (RCO) technique \cite{BCF1, BCF2} to deal with the last term in the RHS of the inequality \eqref{eq:erfinial}. From the NLDE \eqref{eq:Dirac_1D}, we find that $\partial_t \Phi +i\bT \Phi= O(1)$.
This motivates the introduction of the following `twisted variable'
\begin{equation}
    \Psi(t, x)=e^{it\bT}\Phi(t, x), \quad t \geq 0,
\end{equation}
which satisfies the equation $\partial_t \Psi(t,x)=-ie^{it\bT}\Big(V I_2+F(e^{-it\bT}\Psi(t, x))\Big)e^{-it\bT}\Psi(t, x)$.

 We have  
$\|\Psi \|_{L^\infty([0,T];(H_{per}^m)^2)} \lesssim 1$ \  and  $\|\partial_t \Psi \|_{L^\infty([0,T];(H_{per}^m)^2)}   \lesssim 1$ with
\begin{equation}\label{Psi}
    \|\Psi(t_{n+1})-\Psi(t_n)  \|_{H^m}\lesssim \tau,\qquad 0 \leq n \leq \frac{T}{\tau}-1,
\end{equation}
under assumptions (A) and (B).

\noindent
{\textbf{Step 1.}} As introduced in \cite{BCF1}, we choose the cut-off parameter $\tau_0 \in (0, 1)$ and $M_0 = 2\lceil 1/\tau_0 \rceil \in \mathbb{Z}^+$ ($\lceil \cdot \rceil$ is the ceiling function) with $1/\tau_0 \leq M_0/2 < 1+ 1/\tau_0$. Under assumptions (A) and (B), we have the following estimate
\begin{equation}
\left\|P_{M_0}\mathscr{F}\Big(e^{-it_k\bT}P_{M_0}\Psi(t_k)\Big) - \mathscr{F}\Big(e^{-it_k\bT}\Psi(t_k)\Big)\right\|_{H^1}	\lesssim  \tau \tau_0^{m-1}.
\end{equation}
Based on the above estimates,  \eqref{eq:erfinial} implies for $0 \leq n \leq T/\tau-1$,
\begin{equation}
\left\|{\bf e}^{[n+1]}\right\|_{H^1}\lesssim \tau_0^{m-1} + \tau + \tau \sum_{k=0}^{n} \left\|{\bf e}^{[k]}\right\|_{H^1} + \left\| \mathscr{L}^n\right\|_{H^1},
\label{eq:matl}
\end{equation}
where
\begin{equation}
\mathscr{L}^n = \sum_{k=0}^{n} e^{i(k+1)\tau \bT}P_{M_0}\mathscr{F}\left(e^{-it_k\bT}(P_{M_0}\Psi(t_k))\right).
\end{equation}
The quantity $\mathscr{F}(\Phi(t_n))$ admits the decomposition
$$
\mathscr{F}\big(\Phi(t_n)\big)
=
\mathscr{F}_V\big(\Phi(t_n)\big)
+
\mathscr{F}_F\big(\Phi(t_n)\big),
$$
where the components $\mathscr{F}_V\big(\Phi(t_n)\big)$ and
$\mathscr{F}_F\big(\Phi(t_n)\big)$ are defined, respectively, by
\begin{align}
    \mathscr{F}_V\big(\Phi(t_n)\big)&=i\biggl( \int_0^\tau e^{-i(\tau-s)\bT}\ V I_2\ e^{-is\bT}\Phi(t_n)\ ds-\tau  e^{-i\tau \bT} V I_2  \ \Phi(t_n) \biggr),\label{VF}\\
    \mathscr{F}_F\big(\Phi(t_n)\big)&=i\biggl( \int_0^\tau e^{-i(\tau-s)\bT} F(e^{-is\bT }\Phi(t_n)) e^{-is\bT}\Phi(t_n) ds\nonumber\\
        &\qquad -\tau  e^{-i\tau \bT}  F(\Phi(t_n)) \Phi(t_n)\biggr).\label{VF2}
\end{align}
Consequently, $\mathscr{L}^n$ admits the decomposition
$\mathscr{L}^n=\mathscr{L}_V^n+\mathscr{L}_F^n,$
with $\mathscr{L}_V^n$ and $\mathscr{L}_F^n$ defined as follows
\begin{align} \label{L}
   & \mathscr{L}_V^n = \sum_{k=0}^{n} e^{i(k+1)\tau \bT}P_{M_0}\mathscr{F}_V\big(e^{-it_k\bT}(P_{M_0}\Psi(t_k))\big),
    \\&\label{L2}
    \mathscr{L}_F^n = \sum_{k=0}^{n} e^{i(k+1)\tau \bT}P_{M_0}\mathscr{F}_F\big(e^{-it_k\bT}(P_{M_0}\Psi(t_k))\big).
\end{align}

\noindent
{\textbf{Step 2.}} We estimate the two low-frequency contributions  $\mathscr{L}_V^n$ and $\mathscr{L}_F^n$ separately. For $l \in \mathcal{T}_{M_0}$, define the index set $\mathcal{I}_l^{M_0}$ and $\mathcal{J}_l^{M_0}$ associated to $l$ as
\begin{align}
   &\mathcal{I}_l^{M_0} = \left\{(l_1,l_2) \ \vert \ l_1+l_2 = l, (l_1, l_2)  \in\mathcal{T}_{M_0}\times \mathbb{Z}\right\},\\&
    \mathcal{J}_l^{M_0}=\left\{(l_1,l_2,l_3) \ \vert \ -l_1+l_2+l_3 = l,\ l_1, l_2,l_3 \in\mathcal{T}_{M_0}\right\}.
\end{align}
For the subsequent analysis, we define
\begin{equation}
    \Pi_l^+= Q_l\,\text{diag}(1,0) (Q_l)^T,\quad \Pi_l^- = Q_l\,\text{diag}(0,1) (Q_l)^T,
\end{equation}
where $\Pi_l^{\pm}$ denote the projectors onto the eigenspaces of $\Gamma_l$ corresponding to the eigenvalues $\pm\omega_l$, respectively. These projectors satisfy the following identities:
$$
(\Pi_l^\pm)^T=\Pi_l^\pm,\qquad
\Pi_l^++\Pi_l^-=I_2,\qquad
(\Pi_l^\pm)^2=\Pi_l^\pm,\qquad
\Pi_l^{\pm}\Pi_l^{\mp}={\bf 0}.
$$
The spectral decomposition gives
\begin{equation}
e^{-it{\bT}}P_{M_0}\Psi(t_k)=\sum\limits_{l\in\mathcal{T}_{M_0}}\left(e^{-it\omega_l}\Pi_l^++e^{it\omega_l}\Pi_l^-\right)\widehat{\Psi}_l(t_k)e^{i\mu_l(x-a)}.
\end{equation}
In what follows, we analyze the two terms $\mathscr{L}_V^n$ and $\mathscr{L}_F^n$ separately. We begin with the low-Fourier-mode term $\mathscr{L}_V^n$. The definition of $\mathscr{F}_V$ in \eqref{VF} yields the following expansion
\begin{align*}
    &\ e^{i(k+1)\tau \bT}P_{M_0}\left(e^{-i(\tau-s)\bT}Ve^{-is\bT}\big(e^{-it_k\bT}P_{M_0}(\Psi(t_k))\big)\right)
\\= & \sum\limits_{l\in\mathcal{T}_{M_0}}\sum\limits_{(l_1,l_2)\in\mathcal{I}_l^{M_0}}\sum\limits_{\nu_j=\pm, j=1, 2}\mathcal{G}_{k,l,l_1,l_2}^{\nu_1,\nu_2}(s) e^{i\mu_l(x-a)},
    \end{align*}
where  $\mathcal{G}_{k,l,l_1,l_2}^{\nu_1,\nu_2}(s)$ is a function of $s$ defined as
\begin{equation}
\mathcal{G}_{k,l,l_1,l_2}^{\nu_1, \nu_2}(s) =
e^{i(t_k+s)\omega^{\nu_1, \nu_2}_{l,l_1}}\Pi_l^{\nu_1}\widehat{V}_{l_2}\Pi_{l_1}^{\nu_2} \widehat{\Psi}_{l_1}(t_k),
\label{eq:mGdef}
\end{equation}
with $\omega_{l,l_1}^{\nu_1,\nu_2}=\nu_1 \omega_l-\nu_2 \omega_{l_1}$.
Then, the remainder term $\mathscr{L}_V^n$ in \eqref{L} reads
\begin{equation}
\label{eq:reminder-dec}
\mathscr{L}_V^n =i\sum\limits_{k=0}^n
\sum\limits_{l\in\mathcal{T}_{M_0}}\sum\limits_{(l_1,l_2)\in\mathcal{I}_l^{M_0}}\sum\limits_{ \nu_j=\pm, j=1,2} \eta_{k, l, l_1, l_2}^{\nu_1, \nu_2} e^{i\mu_l(x-a)},
\end{equation}
where 
\begin{align}
\eta_{k, l, l_1, l_2}^{\nu_1, \nu_2} & = \int_0^\tau\mathcal{G}_{k,l,l_1,l_2}^{\nu_1, \nu_2}(s)\,ds-\tau\mathcal{G}_{k,l,l_1,l_2}^{\nu_1, \nu_2}(0)  = r_{l,l_1}^{\nu_1, \nu_2}e^{it_k\omega^{\nu_1, \nu_2}_{l,l_1}}c_{k,l,l_1,l_2}^{\nu_1, \nu_2},
\label{eq:Lambdak}
\end{align}
and  the vector coefficients  $c_{k,l,l_1,l_2}^{\nu_1, \nu_2}$  and the scalar coefficient $r_{l,l_1}^{\nu_1, \nu_2}$ are given as
\begin{align}
& c_{k,l,l_1,l_2}^{\nu_1, \nu_2} = 	\Pi_l^{\nu_1}\widehat{V}_{l_2}\Pi_{l_1}^{\nu_2} \widehat{\Psi}_{l_1}(t_k), \\
& r_{l,l_1}^{\nu_1, \nu_2} = \int^{\tau}_0 e^{is\omega^{\nu_1, \nu_2}_{l,l_1}}\;ds-\tau = O(\tau^2 \omega^{\nu_1, \nu_2}_{l,l_1} ).
\label{eq:rest1}
\end{align}
We only need to consider the case $\omega^{\nu_1,\nu_2}_{l, l_1} \neq0$ as $ r_{l, l_1}^{\nu_1, \nu_2} = 0$ if $\omega^{\nu_1, \nu_2}_{l, l_1} = 0$. For $l\in\mathcal{T}_{M_0}$ and $(l_1, l_2)\in\mathcal{I}_l^{M_0}$, we have
\begin{equation}
\vert\omega^{\nu_1,\nu_2}_{l, l_1}\vert \leq 2\omega_{M_0/2}= 2\sqrt{\delta^2+\frac{\mu_{{M_0}/2}^2}{\eps^2}} < 2\sqrt{\delta^2 + \frac{4\pi^2 (1+\tau_0)^2}{\eps^2\tau_0^2(b-a)^2}},
\end{equation}
which implies for $0< \tau \leq \eps\alpha \frac{\pi(b-a)\tau_0}{\sqrt{\delta^2 \eps^2\tau_0^2(b-a)^2 + 4\pi^2(1+\tau_0)^2}}$ ($0 < \tau_0, \alpha < 1$), there holds $\frac{\tau}{2}\vert \omega_{l, l_1}^{\nu_1, \nu_2} \vert \leq \alpha\pi$.

Denoting $S^{\nu_1, \nu_2}_{n, l, l_1}=\sum\limits_{k=0}^n e^{it_k \omega_{l, l_1}^{\nu_1, \nu_2}}$,
for $0< \tau \leq \eps\alpha \frac{\pi(b-a)\tau_0}{\sqrt{\delta^2 \eps^2\tau_0^2(b-a)^2 + 4\pi^2(1+\tau_0)^2}}$, we obtain
\begin{equation}
\label{eq:Sbd1}
 \vert S^{\nu_1, \nu_2}_{n, l, l_1}\vert \leq \frac{1}{\vert \sin(\tau \omega_{l, l_1}^{\nu_1, \nu_2}/2)\vert}\leq\frac{C}{\tau\vert\omega_{l, l_1}^{\nu_1, \nu_2}\vert},\qquad C = \frac{2\alpha\pi}{\sin(\alpha\pi)},\qquad  \forall\   n \ge0.
\end{equation}
An application of the summation-by-parts formula yields
\begin{align}
\sum_{k=0}^n\eta_{k,l,l_1,l_2}^{\nu_1, \nu_2}=
r_{l, l_1}^{\nu_1, \nu_2}\left[\sum_{k=0}^{n-1}S^{\nu_1,\nu_2}_{k,l,l_1} (c_{k,l,l_1,l_2}^{\nu_1, \nu_2} - c_{k+1,l,l_1,l_2}^{\nu_1, \nu_2})  +S_{n, l, l_1}^{\nu_1, \nu_2}c_{n,l,l_1,l_2}^{\nu_1, \nu_2}\right],
\label{eq:lambdasum1}
\end{align}
with
\begin{equation}
c_{k,l,l_1,l_2}^{\nu_1, \nu_2}-c_{k+1,l,l_1,l_2}^{\nu_1, \nu_2} = 	\Pi_l^{\nu_1}\widehat{V}_{l_2}\Pi_{l_1}^{\nu_2}\left(\widehat{\Psi}_{l_1}(t_k) - \widehat{\Psi}_{l_1}(t_{k+1})\right).
\label{eq:cksum1}
\end{equation}
Combining \eqref{eq:rest1}, \eqref{eq:Sbd1}, \eqref{eq:lambdasum1}, and \eqref{eq:cksum1}, we have
\begin{equation}
\left\vert\sum_{k=0}^n \eta_{k,l,l_1,l_2}^{\nu_1, \nu_2}\right\vert
\lesssim   \tau\left\vert\widehat{V}_{l_2}\right\vert\left[\sum\limits_{k=0}^{n-1} \left\vert\widehat{\Psi}_{l_1}(t_k)-\widehat{\Psi}_{l_1}(t_{k+1})\right\vert
 + \left\vert\widehat{\Psi}_{l_1}(t_n)\right\vert\right].\label{eq:sumlambda1}
\end{equation}
Now, we are ready to give the improved estimates for $\mathscr{L}_V^n$. For $l\in\mathcal{T}_{M_0}$ and $(l_1, l_2)\in\mathcal{I}_l^{M_0}$, simple calculations show that $l = l_1 + l_2$ and
\begin{equation}
1+\mu_l^2\lesssim\prod_{j=1}^2\left(1+\mu_{l_j}^2\right).
\label{eq:mlbd1}
\end{equation}
Based on \eqref{eq:reminder-dec}, \eqref{eq:sumlambda1} and \eqref{eq:mlbd1}, we have
\begin{align}
&\ \left\| \mathscr{L}_V^n\right\|^2_{H^1}  \nonumber\\
=& \
\sum\limits_{l\in\mathcal{T}_{M_0}} \left(1+\mu_l^2\right)   \left\vert \sum\limits_{(l_1,l_2)\in\mathcal{I}_l^{M_0}}\sum\limits_{\nu_j = \pm,  j=1, 2}\sum\limits_{k=0}^n\eta_{k, l, l_1, l_2}^{\nu_1,\nu_2}\right\vert^2\nonumber\\
\lesssim&\ \tau^2
\bigg\{\sum_{l\in\mathcal{T}_{M_0}}\bigg(\sum\limits_{(l_1, l_2)\in\mathcal{I}_l^{M_0}}\left\vert \widehat{V}_{l_2}\right\vert \left\vert \widehat{\Psi}_{l_1}(t_n)\right\vert \prod_{j=1}^2\left(1+\mu_{l_j}^2\right) \bigg)^{2}\nonumber\\
& + n \sum\limits_{k=0}^{n-1} \bigg[\sum_{l\in\mathcal{T}_{M_0}}\bigg(\sum\limits_{(l_1, l_2)\in\mathcal{I}_l^{M_0}}
\left\vert\widehat{V}_{l_2}\right\vert\left\vert\widehat{\Psi}_{l_1}(t_k)-\widehat{\Psi}_{l_1}(t_{k+1})\right\vert\prod_{j=1}^2\left(1+\mu_{l_j}^2\right) \bigg)^2\bigg]\bigg\}.\label{eq:sumlambda-21}
\end{align}
In order to estimate each term in the above inequality, we use the auxiliary function $\Theta(x)=\sum_{l\in\mathbb{Z}}(1+\mu_l^2)\left\vert\widehat{\Psi}_l(t_n)\right\vert e^{i\mu_l(x-a)}$, where $\Theta(x)\in H_{\rm per}^{m-2}(\Omega)$ is implied by assumption (B). Moreover, we can prove that $\|\Theta(x)\|_{H^{s}}\lesssim\|\Psi(t_n)\|_{H^{s+2}}$ ($s \leq m-2$). Similarly, we introduce the function $\bU(x) = \sum_{l\in\mathbb{Z}}(1+\mu_l^2)\left\vert \widehat{V}_l\right\vert e^{i\mu_l(x-a)}$, where $\bU(x) \in {W^{m-3, \infty}_{\rm per}} (\Omega)$  can be derived from the assumption (A). We can directly prove that  $\|\bU(x)\|_{W^{0,\infty}}\lesssim\|V(x)\|_{W^{3, \infty}}$. Expanding
\[\bU(x)\Theta(x)=\sum\limits_{l\in\mathbb{Z}}\sum\limits_{l_1+l_2=l} \prod_{j=1}^2 (1+\mu_{l_j}^2) \left\vert\widehat{V}_{l_2}\right\vert\left\vert\widehat{\Psi}_{l_1}(t_n)\right\vert e^{i\mu_l(x-a)},\]
we can obtain
\begin{align}
&\sum_{l\in\mathcal{T}_{M_0}} \bigg(\sum\limits_{(l_1, l_2)\in\mathcal{I}_l^{M_0}}\left\vert\widehat{V}_{l_2}\right\vert\left\vert\widehat{\Psi}_{l_1}(t_n)\right\vert \prod_{j=1}^2\left(1+\mu_{l_j}^2\right) \bigg)^2\nonumber\\
& \leq \ \|\bU(x)\Theta(x)\|^2_{L^2}\lesssim  \| \bU(x)\|_{L^\infty}^2\|\Theta(x)\|_{H^1}^2\lesssim \|V(x)\|_{W^{3, \infty}}^2 \|\Psi(t_n)\|_{H^3}^2\lesssim1.
\label{eq:est_ll}
\end{align}
Using \eqref{Psi},  we  obtain the estimate for $\mathscr{L}_V^n$ similarly as
\begin{align}
\left\| \mathscr{L}_V^n\right\|^2_{H^1} & \ \lesssim \tau^2\|V(x)\|_{W^{3, \infty}}^2 \bigg(\|\Psi(t_n)\|_{H^3}^2 + n\sum\limits_{k=0}^{n-1}
\|\Psi(t_k) - \Psi(t_{k+1})\|_{H^3}^2\bigg)\nonumber \\
&\  \lesssim  \tau^2+ n^2\tau^2(\tau^2)
\lesssim\tau^2.
\label{eq:est-l2l}
\end{align}
Next, we analyze the low-Fourier-mode term $\mathscr{L}_F^n$. In view of the definition of $\mathscr{F}_F$ in \eqref{VF2}, we have the following expansion
\begin{align*}
    &\ e^{i(k+1)\tau \bT}P_{M_0}\bigg(e^{-i(\tau-s)\bT}F\Big(e^{-i(t_k+s)\bT}P_{M_0}\Psi(t_k)\Big)e^{-i(t_k+s)\bT}P_{M_0}\left(\Psi(t_k)\right)\bigg)
\\=&\ \lambda\sum\limits_{l\in\mathcal{T}_{M_0}}\sum\limits_{(l_1,l_2,l_3)\in\mathcal{J}_l^{M_0}}\sum\limits_{\nu_j=\pm, j=0,1,2,3}\mathcal{G}_{k,l,l_1,l_2,l_3}^{\nu_0,\nu_1,\nu_2,\nu_3}(s) e^{i\mu_l(x-a)},
    \end{align*}
where $\mathcal{G}_{k,l,l_1,l_2,l_3}^{\nu_0,\nu_1,\nu_2,\nu_3}(s)$ is a function of $s$ defined as
\begin{equation}
\mathcal{G}_{k,l,l_1,l_2,l_3}^{\nu_0,\nu_1,\nu_2,\nu_3}(s)=
e^{i(t_k+s)\omega^{\nu_0,\nu_1, \nu_2,\nu_3}_{l,l_1,l_2,l_3}}\Pi_l^{\nu_0}\left((\widehat{\Psi}_{l_1}(t_k))^{*}\Pi_{l_1}^{\nu_1}\sigma_3\Pi_{l_2}^{\nu_2}\widehat{\Psi}_{l_2}(t_k)\right)\sigma_3\Pi_{l_3}^{\nu_3}\widehat{\Psi}_{l_3}(t_k),
\label{eq:mGdef2}
\end{equation}
with $\omega_{l,l_1,l_2,l_3}^{\nu_0,\nu_1,\nu_2,\nu_3}=\nu_0\omega_l+\nu_1 \omega_{l_1}-\nu_2 \omega_{l_2}-\nu_3 \omega_{l_3}$.
Then, the remainder term $\mathscr{L}_F^n$ in \eqref{L2} reads
\begin{equation}
\label{eq:reminder-dec2}
\mathscr{L}_F^n =i\lambda\sum\limits_{k=0}^n
\sum\limits_{l\in\mathcal{T}_{M_0}}\sum\limits_{(l_1,l_2,l_3)\in\mathcal{J}_l^{M_0}}\sum\limits_{ \nu_j=\pm, j=0,1,2,3} \eta_{k, l, l_1, l_2,l_3}^{\nu_0,\nu_1, \nu_2,\nu_3} e^{i\mu_l(x-a)},
\end{equation}
where 
\begin{align}
\eta_{k, l, l_1, l_2,l_3}^{\nu_0,\nu_1, \nu_2,\nu_3}  =&\ \int_0^\tau\mathcal{G}_{k,l,l_1,l_2,l_3}^{\nu_0,\nu_1,\nu_2,\nu_3}(s)ds-\tau\mathcal{G}_{k,l,l_1,l_2,l_3}^{\nu_0,\nu_1,\nu_2,\nu_3}(0)\\=&\ r_{l,l_1,l_2,l_3}^{\nu_0,\nu_1, \nu_2,\nu_3}e^{it_k\omega^{\nu_0,\nu_1, \nu_2,\nu_3}_{l,l_1,l_2,l_3}}c_{k,l,l_1,l_2,l_3}^{\nu_0,\nu_1, \nu_2,\nu_3},
\label{eq:Lambdak2}
\end{align}
and  the vector coefficients  $c_{k,l,l_1,l_2,l_3}^{\nu_0,\nu_1, \nu_2,\nu_3}$  and the scalar coefficient $r_{l,l_1,l_2,l_3}^{\nu_0,\nu_1, \nu_2,\nu_3}$ are given as
    \begin{align}\label{eq:rest12}
&c_{k,l,l_1,l_2,l_3}^{\nu_0,\nu_1, \nu_2,\nu_3} = \ 	\Pi_l^{\nu_0}\left((\widehat{\Psi}_{l_1}(t_k))^{*}\Pi_{l_1}^{\nu_1}\sigma_3\Pi_{l_2}^{\nu_2}\widehat{\Psi}_{l_2}(t_k)\right)\sigma_3\Pi_{l_3}^{\nu_3}\widehat{\Psi}_{l_3}(t_k), \\
& r_{l,l_1,l_2,l_3}^{\nu_0,\nu_1, \nu_2,\nu_3} =\  \int^{\tau}_0 e^{is\omega_{l,l_1,l_2,l_3}^{\nu_0,\nu_1,\nu_2,\nu_3}}\;ds-\tau = O(\tau^2 \omega_{l,l_1,l_2,l_3}^{\nu_0,\nu_1,\nu_2,\nu_3}).
\end{align}
We only need to consider the case
$\omega^{\nu_0,\nu_1,\nu_2,\nu_3}_{l, l_1,l_2,l_3} \neq0$ as  $ r_{l, l_1,l_2,l_3}^{\nu_0,\nu_1, \nu_2,\nu_3} = 0$  if $\omega^{\nu_0,\nu_1, \nu_2,\nu_3}_{l, l_1,l_2,l_3} = 0$. For $l\in\mathcal{T}_{M_0}$ and $(l_1, l_2,l_3)\in\mathcal{J}_l^{M_0}$, we have
\begin{equation}
\vert\omega^{\nu_0,\nu_1,\nu_2,\nu_3}_{l, l_1,l_2,l_3}\vert \leq 4\omega_{M_0/2}= 4\sqrt{\delta^2+\frac{\mu_{{M_0}/2}^2}{\eps^2}} < 4\sqrt{\delta^2 + \frac{4\pi^2 (1+\tau_0)^2}{\eps^2\tau_0^2(b-a)^2}},
\end{equation}
which implies $\frac{\tau}{2}\vert \omega_{l, l_1,l_2,l_3}^{\nu_0,\nu_1, \nu_2,\nu_3} \vert \leq \alpha\pi$ when $0< \tau \leq \eps\alpha \frac{\pi(b-a)\tau_0}{2  \sqrt{\delta^2 \eps^2\tau_0^2(b-a)^2 + 4\pi^2(1+\tau_0)^2}}$ ( $0 < \tau_0, \alpha < 1$). For $0< \tau \leq \eps\alpha \frac{\pi(b-a)\tau_0}{2 \sqrt{\delta^2 \eps^2\tau_0^2(b-a)^2 + 4\pi^2(1+\tau_0)^2}}$, we obtain the following
\begin{equation}
\label{eq:Sbd12}
 \vert S^{\nu_0,\nu_1, \nu_2,\nu_3}_{n, l, l_1,l_2,l_3}\vert \leq \frac{1}{\vert \sin(\tau \omega_{l, l_1,l_2,l_3}^{\nu_0,\nu_1, \nu_2,\nu_3}/2)\vert}\leq\frac{C}{\tau\vert\omega_{l, l_1,l_2,l_3}^{\nu_0,\nu_1, \nu_2,\nu_3}\vert},\quad C = \frac{2\alpha\pi}{\sin(\alpha\pi)},\quad  \forall \ n\ge0,
\end{equation}
where $S^{\nu_0,\nu_1, \nu_2,\nu_3}_{n, l, l_1,l_2,l_3}=\sum\limits_{k=0}^n e^{it_k \omega_{l, l_1,l_2,l_3}^{\nu_0,\nu_1, \nu_2,\nu_3}}$. Using the summation-by-parts, we find that
\begin{equation}
  \begin{aligned}[b]
&\sum_{k=0}^n\eta_{k,l,l_1,l_2,l_3}^{\nu_0,\nu_1, \nu_2,\nu_3}\\= &\ 
r_{l, l_1,l_2,l_3}^{\nu_0,\nu_1, \nu_2,\nu_3}\left[\sum_{k=0}^{n-1}S^{\nu_0,\nu_1,\nu_2,\nu_3}_{k,l,l_1,l_2,l_3} (c_{k,l,l_1,l_2,l_3}^{\nu_0,\nu_1, \nu_2,\nu_3} - c_{k+1,l,l_1,l_2,l_3}^{\nu_0,\nu_1, \nu_2,\nu_3})  +S_{n, l, l_1,l_2,l_3}^{\nu_0,\nu_1, \nu_2,\nu_3}c_{n,l,l_1,l_2,l_3}^{\nu_0,\nu_1, \nu_2,\nu_3}\right]
\label{eq:lambdasum12}
\end{aligned}  
\end{equation}
with
\begin{equation}
\begin{aligned}[b]
&\ c_{k,l,l_1,l_2,l_3}^{\nu_0,\nu_1, \nu_2,\nu_3}-c_{k+1,l,l_1,l_2,l_3}^{\nu_0,\nu_1, \nu_2,\nu_3}\\
=&\ 	\Pi_l^{\nu_0}\left(\big(\widehat{\Psi}_{l_1}(t_k)-\widehat{\Psi}_{l_1}(t_{k+1})\big)^{*}\Pi_{l_1}^{\nu_1}\sigma_3\Pi_{l_2}^{\nu_2}\widehat{\Psi}_{l_2}(t_k)\right)\sigma_3\Pi_{l_3}^{\nu_3}\widehat{\Psi}_{l_3}(t_k)
\\&\ +	\Pi_l^{\nu_0}\left(\big(\widehat{\Psi}_{l_1}(t_{k+1})\big)^{*}\Pi_{l_1}^{\nu_1}\sigma_3\Pi_{l_2}^{\nu_2}\big(\widehat{\Psi}_{l_2}(t_k)-\widehat{\Psi}_{l_2}(t_{k+1})\big)\right)\sigma_3\Pi_{l_3}^{\nu_3}\widehat{\Psi}_{l_3}(t_k)
\\&\ + 	\Pi_l^{\nu_0}\left(\big(\widehat{\Psi}_{l_1}(t_{k+1})\big)^{*}\Pi_{l_1}^{\nu_1}\sigma_3\Pi_{l_2}^{\nu_2}\widehat{\Psi}_{l_2}(t_{k+1})\right)\sigma_3\Pi_{l_3}^{\nu_3}\big(\widehat{\Psi}_{l_3}(t_k)-\widehat{\Psi}_{l_3}(t_{k+1})\big)
\end{aligned}
\label{eq:cksum12}
\end{equation}
Combining \eqref{eq:rest12}, \eqref{eq:Sbd12}, \eqref{eq:lambdasum12}, and \eqref{eq:cksum12}, we have
\begin{equation}
\begin{aligned}[b]
\left\vert\sum_{k=0}^n \eta_{k,l,l_1,l_2,l_3}^{\nu_0,\nu_1, \nu_2,\nu_3}\right\vert
\lesssim    
&\ \tau\sum\limits_{k=0}^{n-1}
    \bigg(
\vert \widehat{\Psi}_{l_1}(t_k)-\widehat{\Psi}_{l_1}(t_{k+1})\vert\vert\widehat{\Psi}_{l_2}(t_k)\vert\vert\widehat{\Psi}_{l_3}(t_k)\vert
\\&\ +\vert \widehat{\Psi}_{l_1}(t_{k+1})\vert\vert\widehat{\Psi}_{l_2}(t_k)-\widehat{\Psi}_{l_2}(t_{k+1})\vert\vert\widehat{\Psi}_{l_3}(t_k)\vert
\\&\ +\vert \widehat{\Psi}_{l_1}(t_{k+1})\vert\vert\widehat{\Psi}_{l_2}(t_{k+1})\vert\vert\widehat{\Psi}_{l_3}(t_k)-\widehat{\Psi}_{l_3}(t_{k+1})\vert\bigg) 
\\&\ +\tau\vert \widehat{\Psi}_{l_1}(t_n)\vert\vert\widehat{\Psi}_{l_2}(t_n)\vert\vert\widehat{\Psi}_{l_3}(t_n)\vert.
\end{aligned}
\label{eq:sumlambda12}
\end{equation}
Now, we are ready to give the improved estimates for $\mathscr{L}_F^n$. For $l\in\mathcal{T}_{M_0}$ and $(l_1, l_2,l_3)\in\mathcal{J}_l^{M_0}$, simple calculations show ( $l = -l_1 + l_2+l_3$ )
\begin{equation}
1+\mu_l^2\lesssim\prod_{j=1}^3\left(1+\mu_{l_j}^2\right).
\label{eq:mlbd12}
\end{equation}
Based on \eqref{eq:reminder-dec2}, \eqref{eq:sumlambda12} and \eqref{eq:mlbd12}, we have
\begin{align}
&\left\| \mathscr{L}_F^n\right\|^2_{H^1}\nonumber\\
&= \
\lambda^2\sum\limits_{l\in\mathcal{T}_{M_0}} \left(1+\mu_l^2\right)   \big\vert \sum\limits_{(l_1,l_2,l_3)\in\mathcal{J}_l^{M_0}}\sum\limits_{\nu_j = \pm,  j=0,1, 2,3}\sum\limits_{k=0}^n\eta_{k, l, l_1, l_2,l_3}^{\nu_0,\nu_1,\nu_2,\nu_3}\big\vert^2\nonumber\\
&\ \lesssim  \tau^2
\bigg\{\sum_{l\in\mathcal{T}_{M_0}}\bigg(\sum\limits_{(l_1, l_2,l_3)\in\mathcal{J}_l^{M_0}}\vert \widehat{\Psi}_{l_1}(t_n)\vert\vert\widehat{\Psi}_{l_2}(t_n)\vert\vert\widehat{\Psi}_{l_3}(t_n)\vert \prod_{j=1}^3\left(1+\mu_{l_j}^2\right) \bigg)^{2}\nonumber\\
&\ + n \sum\limits_{k=0}^{n-1} \bigg[\sum_{l\in\mathcal{T}_{M_0}}\Bigg(\bigg(\sum\limits_{(l_1, l_2,l_3)\in\mathcal{J}_l^{M_0}}\nonumber\\
&\qquad\vert \widehat{\Psi}_{l_1}(t_k)-\widehat{\Psi}_{l_1}(t_{k+1})\vert\vert\widehat{\Psi}_{l_2}(t_k)\vert\vert\widehat{\Psi}_{l_3}(t_k)\vert\prod_{j=1}^3\left(1+\mu_{l_j}^2\right) \bigg)^2\nonumber\\
&\ +\bigg(\sum\limits_{(l_1, l_2,l_3)\in\mathcal{J}_l^{M_0}}
\vert \widehat{\Psi}_{l_1}(t_{k+1})\vert\vert\widehat{\Psi}_{l_2}(t_k)-\widehat{\Psi}_{l_2}(t_{k+1})\vert\vert\widehat{\Psi}_{l_3}(t_k)\vert\prod_{j=1}^3\left(1+\mu_{l_j}^2\right) \bigg)^2\nonumber\\
&\ +\bigg(\sum\limits_{(l_1, l_2,l_3)\in\mathcal{J}_l^{M_0}}
\vert \widehat{\Psi}_{l_1}(t_{k+1})\vert\vert\widehat{\Psi}_{l_2}(t_{k+1})\vert\vert\widehat{\Psi}_{l_3}(t_k)-\widehat{\Psi}_{l_3}(t_{k+1})\vert\nonumber\\
&\qquad\prod_{j=1}^3\left(1+\mu_{l_j}^2\right) \bigg)^2\Bigg)
\bigg]\bigg\}.\label{eq:sumlambda-212}
\end{align}
We use the auxiliary function $\Theta(x)=\sum_{l\in\mathbb{Z}}(1+\mu_l^2)\left\vert\widehat{\Psi}_l(t_n)\right\vert e^{i\mu_l(x-a)}$, where $\Theta(x)\in H_{\rm per}^{m-2}(\Omega)$ implied by assumption (B) and  $\|\Theta(x)\|_{H^{s}} \lesssim\|\Psi(t_n)\|_{H^{s+2}}$ ($s \leq m-2$). Expanding
$$\vert\Theta(x)\vert^2\Theta(x)=\sum\limits_{l\in\mathbb{Z}}\sum\limits_{-l_1+l_2+l_3=l} \prod_{j=1}^3 \left(\left(1+\mu_{l_j}^2\right)\vert\widehat{\Psi}_{l_j}(t_n)\vert\right) e^{i\mu_l(x-a)},$$
we obtain
\begin{align}
&\sum_{l\in\mathcal{T}_{M_0}} \bigg(\sum\limits_{(l_1, l_2,l_3)\in\mathcal{J}_l^{M_0}}\vert\widehat{\Psi}_{l_1}(t_n)\vert \vert\widehat{\Psi}_{l_2}(t_n)\vert \vert\widehat{\Psi}_{l_3}(t_n)\vert \prod_{j=1}^3\left(1+\mu_{l_j}^2\right) \bigg)^2\nn\\
\leq &\ \sum_{l\in\mathbb{Z}} \bigg(\sum\limits_{-l_1+ l_2+l_3=l,l_j\in\mathbb{Z}}\vert\widehat{\Psi}_{l_1}(t_n)\vert \vert\widehat{\Psi}_{l_2}(t_n)\vert \vert\widehat{\Psi}_{l_3}(t_n)\vert \prod_{j=1}^3\left(1+\mu_{l_j}^2\right) \bigg)^2\nn\\
 = &\ \|\vert\Theta(x)\vert^2\Theta(x)\|^2_{L^2}\lesssim \|\Theta(x)\|_{H^1}^6\lesssim \|\Psi(t_n)\|_{H^3}^6\lesssim1.
\label{eq:est_ll2}
\end{align}
Using \eqref{Psi}, we obtain the estimate for each term in $\mathscr{L}_F^n$ similarly as
\begin{align}
& \left\| \mathscr{L}_F^n\right\|^2_{H^1}\nn \\ \lesssim &  \        \tau^2 \bigg[\|\Psi(t_n)\|_{H^3}^6 + n\sum\limits_{k=0}^{n-1}
\|\Psi(t_k) - \Psi(t_{k+1})\|_{H^3}^2(\|\Psi(t_k)\|_{H^3}+\|\Psi(t_{k+1})\|_{H^3})^4\bigg]\nn \\
 \lesssim &\    \tau^2+ n^2\tau^2(\tau^2)
\lesssim\tau^2.\label{eq:est-l2l2}
\end{align}
Combining the results from \eqref{eq:matl}, \eqref{eq:est-l2l} and  \eqref{eq:est-l2l2}, we have
\begin{equation}
    \left\|\mathbf{e}^{[n+1]}\right\|_{H^1}\lesssim\tau_0^{m-1}+\tau+\tau \sum_{k=0}^n\left\|\mathbf{e}^{[k]}\right\|_{H^1},\quad 0\leq n \leq \frac{T}{\tau}-1.
\end{equation}
The discrete Gronwall inequality yields
\begin{equation}
    \left\|\mathbf{e}^{[n+1]}\right\|_{H^1}\lesssim \tau_0^{m-1}+\tau,\quad 0\leq n\leq \frac{T}{\tau}-1.
\end{equation}
This proves \eqref{eq:semi-cutoff} in Theorem~\ref{thm:S2}.

	\section{Numerical Results}\label{sec:num}
In this section, we present some numerical results of the LTFP method to solve the NLDE \eqref{eq:Dirac_21} using various spatial mesh sizes $h$, the time step $\tau$, and the parameter $0<\eps\leq 1$ to illustrate the improved error bounds established above.

\subsection{Numerical results in 1D}
In this subsection, we conduct numerical experiments for the NLDE \eqref{eq:Dirac_21} in the massless nonrelativistic regime, taking $\delta=1$, $\lambda=1$, and $V(x)=2\sin(x)$.

To assess the numerical accuracy of the method, we prescribe the following initial data
\begin{equation}
    \Phi_0(x)=\left(\frac{\sin(x)}{2+\cos(x)^{2}},\,\frac{\cos(x)}{2+\sin(x)^2}\right)^{T},
\qquad x\in (-\pi,\pi ).
\end{equation}
Since the exact solution is unavailable, we compute a numerical reference solution using the LTFP method with a sufficiently small time step $\tau_e=2^{-18}$ and a fine spatial mesh size $h_e=2\pi/2^8$. Let $\Phi^n=(\Phi_0^n,\Phi_1^n,...,\Phi_M^n)^T$ be the numerical solution obtained by the LTFP method with a given mesh size $h$, the time-step size $\tau$ and the parameter $\eps$\ at time $t=t_n$. To evaluate the convergence, we define the numerical error as
\begin{equation*}
    \begin{aligned}
        e(t_n)&=\|\Phi^n-\Phi(t_n,\cdot )\|_{H^1}\\
              &=\sqrt{\,h\sum_{j=0}^{M-1}\bigl|\Phi^{n}_j-\Phi(t_n,x_j)\bigr|^{2}+h\sum_{j=0}^{M-1}\bigl|(\Phi')^{n}_j-\Phi'(t_n,x_j)\bigr|^{2}},
    \end{aligned}
\end{equation*}
where $$ (\Phi')^{n}_j= i\sum_{l\in\mathcal{T}_M}\mu_l\widetilde{(\Phi^n)}_le^{i\mu_l(x_j-a)},j=0,1,...,M-1.$$

\begin{figure}[htbp]
	\centering
	\includegraphics[width=0.45\linewidth]{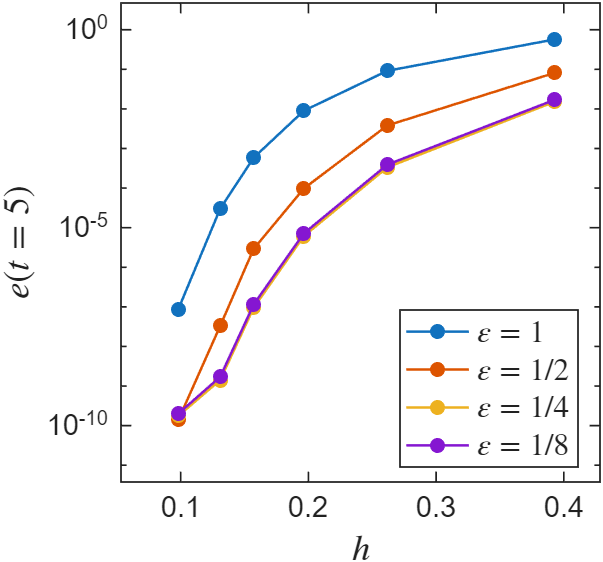}
	\hfill
	\includegraphics[width=0.45\linewidth]{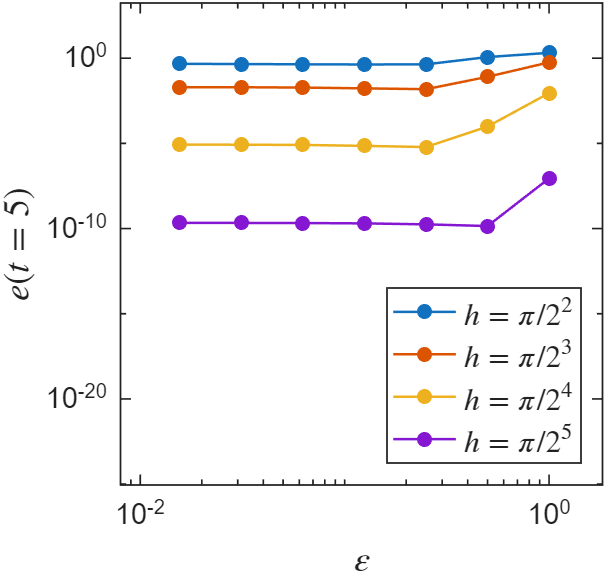}
	\caption{  Spatial errors of the LTFP method for the NLDE \eqref{eq:Dirac_21} in 1D at $t=5$.}
    \label{FIG.1}
\end{figure}

\begin{figure}[htbp]
	\centering
	\includegraphics[width=0.45\linewidth]{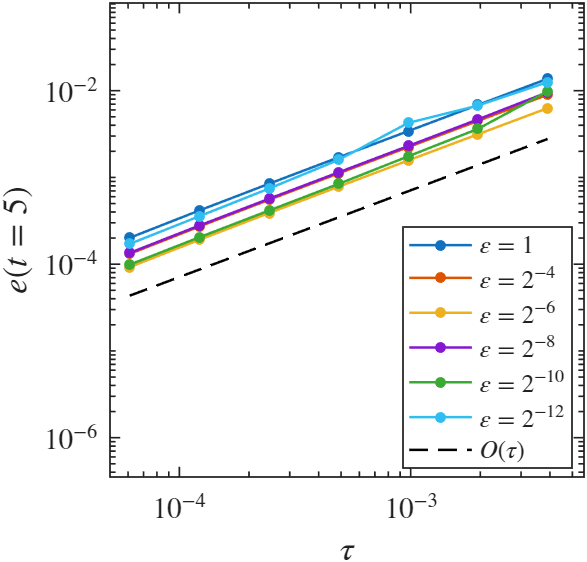}
	\hfill
	\includegraphics[width=0.45\linewidth]{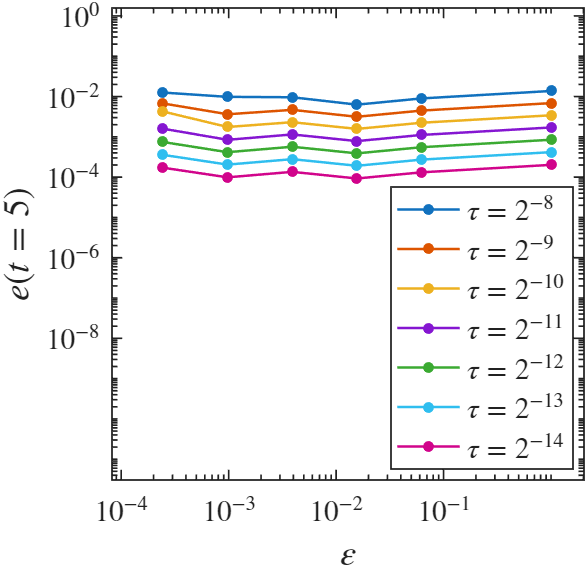}
	\caption{  Temporal errors of the LTFP method for the NLDE \eqref{eq:Dirac_21} in 1D at $t=5$.}	
    \label{FIG.2}
\end{figure}

Figures~\ref{FIG.1} and~\ref{FIG.2} present, respectively, the spatial and temporal errors of the LTFP method applied to the one-dimensional NLDE \eqref{eq:Dirac_21} at the final time $t_n=T=5$. Figure~\ref{FIG.1} demonstrates that the LTFP method achieves spectral accuracy in space and that the spatial discretization errors are independent of the parameter $\eps$.  From Fig. \ref{FIG.2}(a), we observe the first-order convergence of the LTFP method in time for each fixed $\eps$. Figure \ref{FIG.2}(b) again validates  the improved error bounds at $O(\tau)$ on a fixed time interval.

\subsection{Numerical results in 2D}
In this subsection, we conduct numerical experiments for the nonlinear Dirac equation \eqref{eq:Dirac_21} in the massless nonrelativistic regime, taking $\delta=1$, $\lambda=1$, and $V(\bx)=2\bigl(\sin(x)+\sin(y)\bigr)$.

\begin{figure}[htbp]
	\centering
	\includegraphics[width=0.43\linewidth]{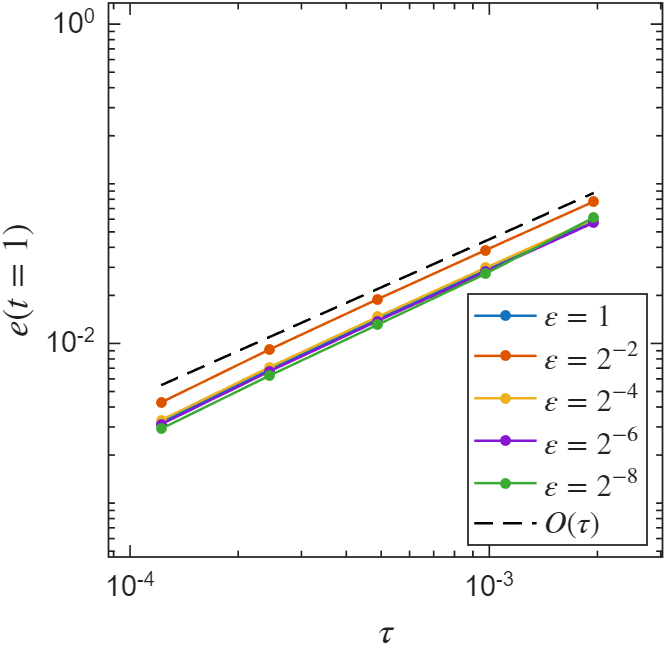}
	\hfill
	\includegraphics[width=0.43\linewidth]{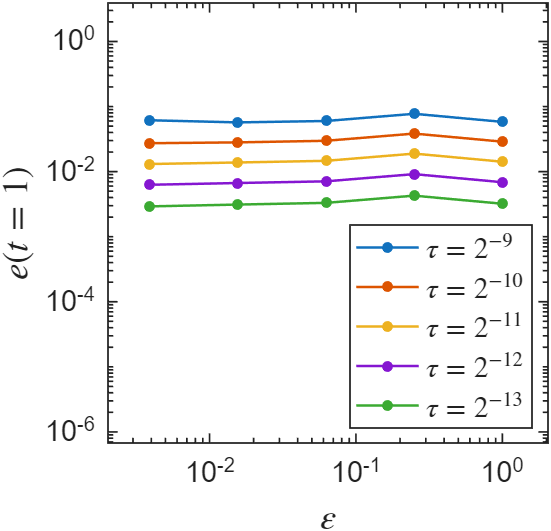}
	\caption{  Temporal errors of the LTFP method for the NLDE \eqref{eq:Dirac_21} in 2D at $t=1$.}
	\label{FIG.3}
\end{figure}

To assess the numerical accuracy of the method, we prescribe the following initial data
\begin{equation}
    \Phi_0(x,y)=\left(\sin(x)+\cos(y),\,\frac{1}{1+\cos^2(x)}+\cos(y)\right)^{T},
\qquad    (x,y)\in(-\pi, \pi)^2.
\end{equation}
Since the exact solution is unavailable, we compute a numerical reference solution using the LTFP method with a sufficiently small time step $\tau_e=2^{-16}$ and a fine spatial mesh size $h_e=2\pi/2^8$. Figure~\ref{FIG.3} presents the temporal errors of the LTFP method applied to the two-dimensional NLDE \eqref{eq:Dirac_21} at the final time $t_n=T=1$. For each fixed $\eps$, Figure~\ref{FIG.3}(a) demonstrates the first-order temporal convergence of the LTFP method, while Figure~\ref{FIG.3}(b) further confirms the improved $O(\tau)$ error bound on an $O(1)$ time interval.

\subsection{Dynamics of the nonlinear Dirac equation in 2D}
In this subsection, we employ the LTFP method to investigate the dynamics of the two-dimensional NLDE \eqref{eq:Dirac_21} subject to a honeycomb lattice potential. Specifically, we consider the following time-independent honeycomb lattice potential:
\begin{equation}
\begin{aligned}[b]
V(\bx)
={}&
\cos\left(\frac{4\pi}{\sqrt{3}}\mathbf{e}_1\cdot\bx\right)
+\cos\left(\frac{4\pi}{\sqrt{3}}\mathbf{e}_2\cdot\bx\right) 
+\cos\left(\frac{4\pi}{\sqrt{3}}\mathbf{e}_3\cdot\bx\right),
\end{aligned}
\end{equation}
where
\begin{equation}
\mathbf{e}_1=(-1,0)^T,\qquad
\mathbf{e}_2=\left(\frac{1}{2},\frac{\sqrt{3}}{2}\right)^T,\qquad
\mathbf{e}_3=\left(\frac{1}{2},-\frac{\sqrt{3}}{2}\right)^T.
\end{equation}

\begin{figure}[!htb]
    \centering

    \begin{minipage}[b]{0.31\textwidth}
        \centering
        \includegraphics[width=\linewidth]{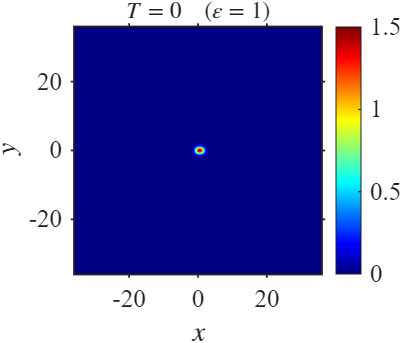}      
    \end{minipage}
    \hfill
    \begin{minipage}[b]{0.31\textwidth}
        \centering
        \includegraphics[width=\linewidth]{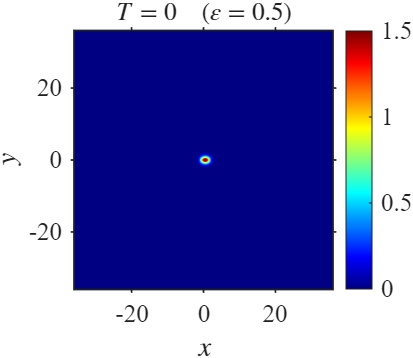}       
    \end{minipage}
    \hfill
    \begin{minipage}[b]{0.31\textwidth}
        \centering
        \includegraphics[width=\linewidth]{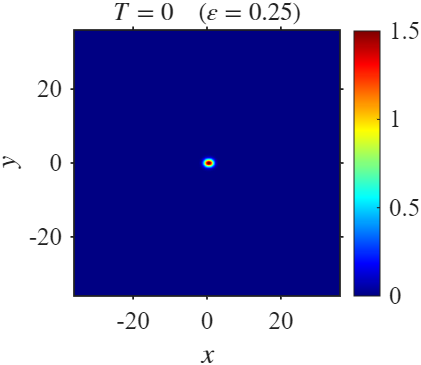}      
    \end{minipage}

    \vspace{0.08cm}
    
    \begin{minipage}[b]{0.31\textwidth}
        \centering
        \includegraphics[width=\linewidth]{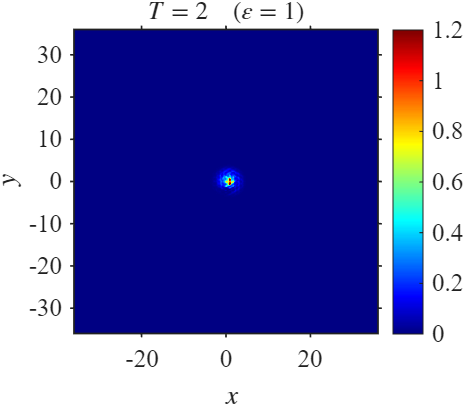}      
    \end{minipage}
    \hfill
    \begin{minipage}[b]{0.31\textwidth}
        \centering
        \includegraphics[width=\linewidth]{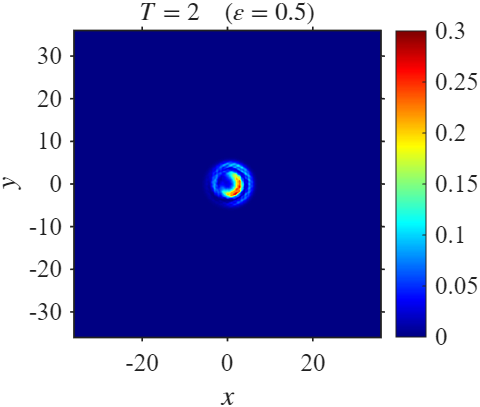}       
    \end{minipage}
    \hfill
    \begin{minipage}[b]{0.31\textwidth}
        \centering
        \includegraphics[width=\linewidth]{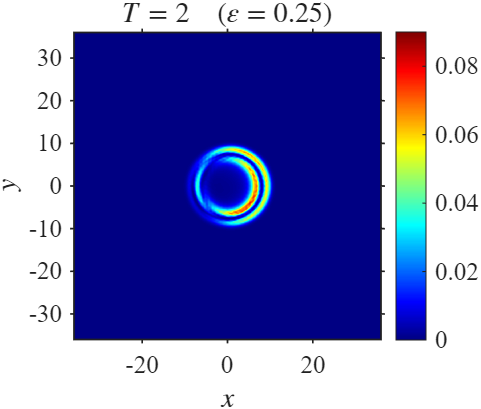}      
    \end{minipage}

    \vspace{0.08cm}
    \begin{minipage}[b]{0.31\textwidth}
        \centering
        \includegraphics[width=\linewidth]{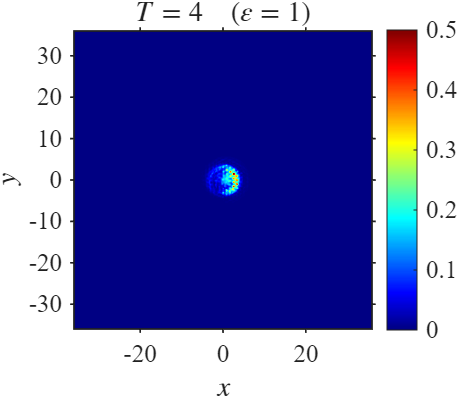}      
    \end{minipage}
    \hfill
    \begin{minipage}[b]{0.31\textwidth}
        \centering
        \includegraphics[width=\linewidth]{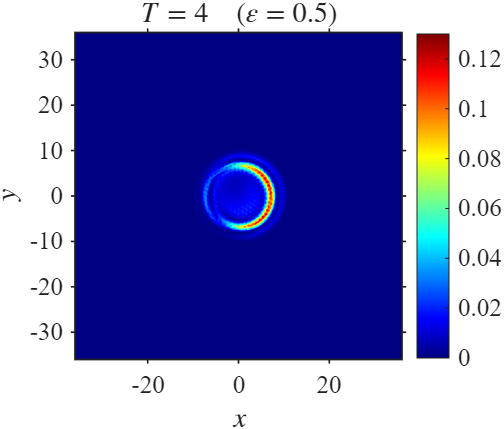}       
    \end{minipage}
    \hfill
    \begin{minipage}[b]{0.31\textwidth}
        \centering
        \includegraphics[width=\linewidth]{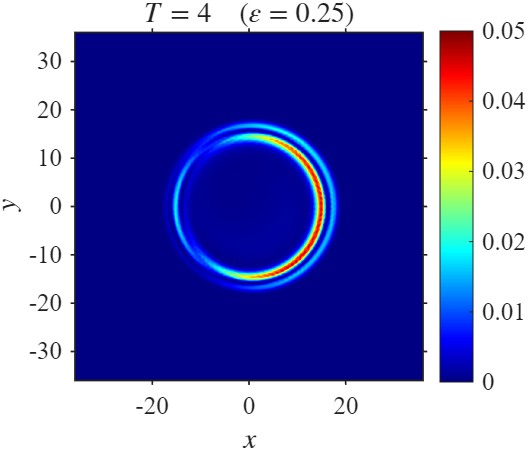}      
    \end{minipage}

    \vspace{0.08cm}

    \begin{minipage}[b]{0.31\textwidth}
        \centering
        \includegraphics[width=\linewidth]{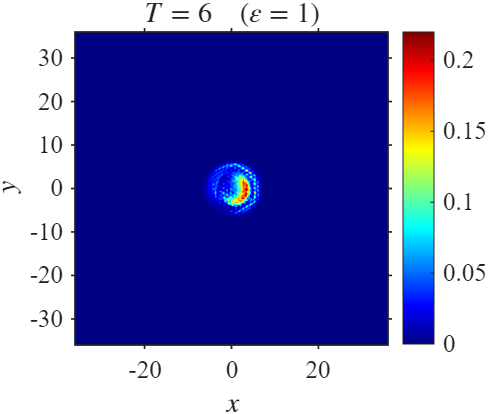}      
    \end{minipage}
    \hfill
    \begin{minipage}[b]{0.31\textwidth}
        \centering
        \includegraphics[width=\linewidth]{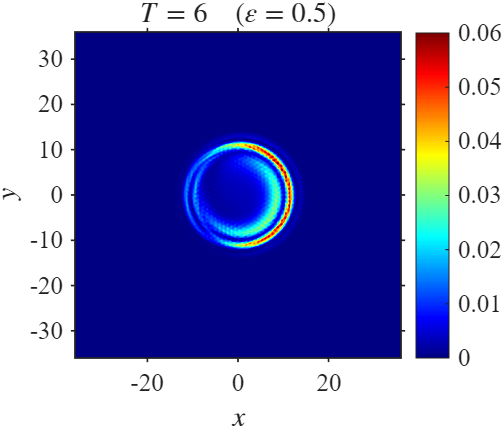}       
    \end{minipage}
    \hfill
    \begin{minipage}[b]{0.31\textwidth}
        \centering
        \includegraphics[width=\linewidth]{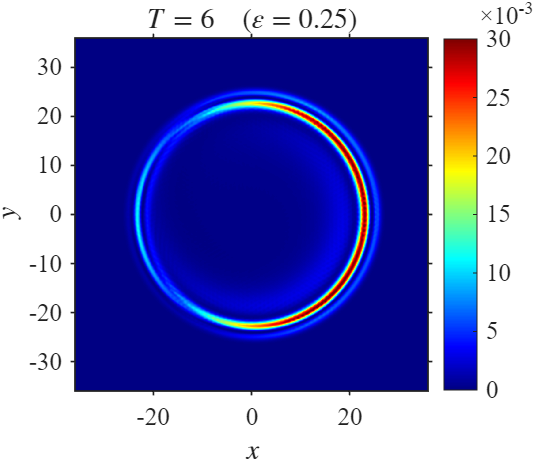}      
    \end{minipage}

    \vspace{0.08cm}
    \begin{minipage}[b]{0.31\textwidth}
        \centering
        \includegraphics[width=\linewidth]{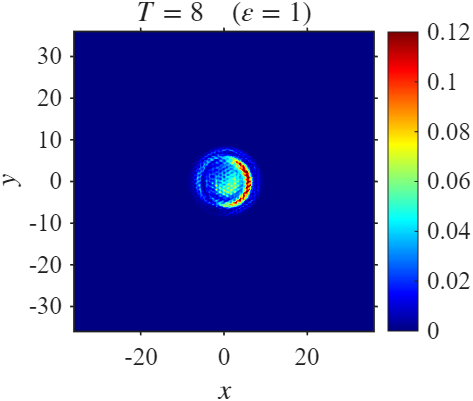}      
    \end{minipage}
    \hfill
    \begin{minipage}[b]{0.31\textwidth}
        \centering
        \includegraphics[width=\linewidth]{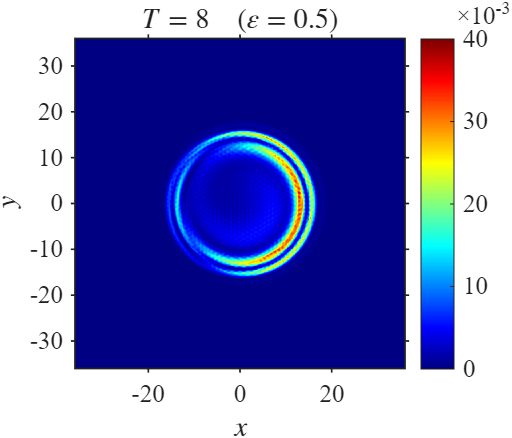}       
    \end{minipage}
    \hfill
    \begin{minipage}[b]{0.31\textwidth}
        \centering
        \includegraphics[width=\linewidth]{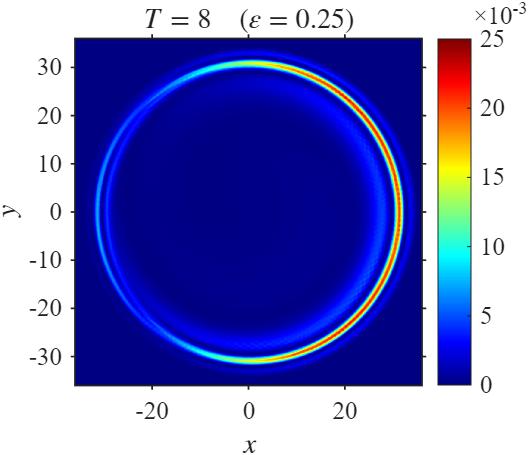}      
    \end{minipage}
   \caption{Dynamics of the total density $\rho$ of the NLDE in 2D with different time $T$ and different $\eps$: $\eps=1$ (left column), $\eps=0.5$ (middle column), and $\eps=0.25$ (right column).}
    \label{fig:rho_dynamics}
\end{figure}

The initial data are chosen as 
\begin{equation}
    \Phi_0(x,y)=\left(e^{-\frac{x^2+y^2}{2}}, e^{-\frac{(x-1)^2+y^2}{2}}\right)^T,  \qquad    (x,y)\in (-36,36)^2.
\end{equation}
We set $\delta=1$ and $\lambda=1$ and solve the problem numerically using the LTFP method with the spatial mesh size $h=1/16$ and the time step $\tau=0.01$. Figure~\ref{fig:rho_dynamics} plots the total density $\rho(t,\bx)$ associated with the solution of the NLDE \eqref{eq:Dirac_21} at several time instances and for different values of $\eps$.
From Figure~\ref{fig:rho_dynamics}, we observe that the dynamics of the nonlinear Dirac equation with the honeycomb lattice potential depends strongly on $\eps$. For $\eps=1$, the density remains relatively localized around the initial wave packet and exhibits irregular local oscillatory structures (cf. left column in Figure~\ref{fig:rho_dynamics}). These local fluctuations are mainly associated with the coupling between the two spinor components, the honeycomb lattice modulation, and the density-imbalance-dependent nonlinear self-interaction. As $\eps$ decreases, the density spreads over a larger spatial region and forms a clearer ring-shaped wave front. In addition, when $\eps$ becomes half, the wave speed doubles (cf. each row in Figure~\ref{fig:rho_dynamics} for different $\eps$), which again confirms the solution propagates wave with wave speed at $O(1/\eps)$.

Moreover, the maximal density decreases as $\eps$ becomes smaller (cf. each row in Figure~\ref{fig:rho_dynamics}), indicating that the wave packet is redistributed over a wider region during propagation. Although nonlinear self-interaction introduces additional modulation, the dominant $\eps$-dependent propagation behavior remains evident. These results show that the Lie-Trotter splitting method with Fourier pseudospectral discretization captures the main dynamics of the two-dimensional nonlinear Dirac equation accurately.
    \FloatBarrier
    
	\section{Conclusions} \label{sec:conclude}
	In this work, we rigorously established improved error bounds for the Lie-Trotter splitting method applied to the nonlinear Dirac equation in the massless nonrelativistic regime. By employing the regularity compensation oscillation (RCO) technique, we established  improved error bounds for semi-discretization and full-discretization of order $O(\tau)$ and $O(h^{m-1}+\tau)$, respectively, especially when $0<\eps\ll 1$. The numerical results in 1D and 2D agreed well with the theoretical results and suggested that the error bounds are sharp. In addition, the Lie-Trotter splitting method was used to study the dynamics of the nonlinear Dirac equation in the massless nonrelativistic regime with honeycomb lattice potential in 2D and some interesting phenomena were observed.

	\section*{Acknowledgment}
	The authors acknowledge the support from the National Key R$\&$D Program of China (No. 2024YFA1016100), Tianyuan Fund for Mathematics of the National Natural Science Foundation of China (No. 12326102), National Natural Science Foundation of China (No. 12401539).

\end{document}